\documentclass[11pt]{article}
\usepackage{amsfonts,amssymb,amsmath,pdfsync,bbm,color}
\hsize \textwidth
\advance \hsize by -\marginparwidth
\evensidemargin \oddsidemargin
\usepackage{amssymb}
\advance\hoffset by 5mm
\newcommand{\pdr}[2]{\frac{\partial{#1}}{\partial{#2}}}
\newcommand{\Rm}{{\mathbb R}}
\newcommand{\eps}{\varepsilon}
\newcommand{\commentout}[1]{}

\newcommand{\be}{\begin{equation}}
\newcommand{\ee}{\end{equation}}
\newcommand{\one}{{\mathbbm 1}}

\newtheorem{thm}{Theorem}[section]
\newtheorem{lem}[thm]{Lemma}
\newtheorem{cor}[thm]{Corollary}
\newtheorem{prop}[thm]{Proposition}

\newtheorem{rem}[thm]{Remark}
\numberwithin{equation}{section} 

\newcommand{\pdrr}[2]{\frac{\partial^2{#1}}{\partial{#2}^2}}

\newcommand{\farc}{\frac}

\newcommand{\bal}{\begin{aligned}}
\newcommand{\enbal}{\end{aligned}}

\newcommand{\Pm}{{\mathbb P}}
\newcommand{\argmax}{\hbox{argmax}}
\newcommand{\tw}{\tilde w}
\newcommand{\etal}{\eta_\ell}
\newcommand{\el}{e_\ell}
\newcommand{\cG}{{\cal G}}
\newcommand{\bF}{{\bar F}}
\newcommand{\tel}{\tilde e_\ell}

\newcommand{\bcR}{\bar{\cal R}}
\newcommand{\un}{\underline}

\newcommand{\bG}{{\bar G}}

\newcommand{\vphi}{\varphi}
\newcommand{\ifnty}{\infty}

\usepackage[colorlinks=true,pdfpagemode=UseNone,urlcolor=blue,linkcolor=blue,citecolor=blue,breaklinks=true]{hyperref}

\begin{document}

\title{Diffusion of knowledge and the  lottery society}

\author{Henri Berestycki\footnote{(a) EHESS, CAMS, 54 boulevard Raspail, F-75006 Paris, France, 
(b) Department of Mathematics, University of Maryland, College Park, USA, and (c) Institute of Advanced Study, 
Hong Kong University of Science and Technology 
Hong Kong; hb@ehess.fr}  \and Alexei Novikov\footnote{Department of Mathematics, The Pennsylvania State University, 
University Park, PA 16802, USA; novikov@psu.edu} 
\and Jean-Michel Roquejoffre\footnote{Institut de Math\'ematiques (UMR CNRS 5219), Universit\'e Paul Sabatier, 118 Route de Narbonne, 31062 Toulouse Cedex, France; jean-michel.roquejoffre@math.univ-toulouse.fr} \and 
Lenya Ryzhik\footnote{Department of Mathematics, Stanford University, Stanford, CA 94305, USA; ryzhik@stanford.edu}}

\maketitle

\begin{abstract}
The Lucas-Moll system is a mean-field game type model describing the growth of an economy by means of
 diffusion of knowledge.
The individual agents in the economy advance their knowledge by learning from each other and via
internal innovation.
Their cumulative distribution function
satisfies a forward in time nonlinear non-local reaction-diffusion type equation. On the other hand,
the learning strategy of the agents is based
on the solution to a backward in time nonlocal Hamilton-Jacobi-Bellman equation
that is coupled to the aforementioned equation for the agents density. Together, these equations
form a system of the mean-field game type. 
When the learning rate is sufficiently large, existence 
of balanced growth path solutions to the Lucas-Moll system was proved in~\cite{PRV,Porretta-Rossi}.
Here, we analyze a complementary regime where the balanced growth paths do not exist. The main result is 
a long time  convergence theorem.  
Namely, the solution to the initial-terminal value problem behaves in such a way that 
at large times 
an overwhelming majority of the  agents spend no time producing at all
and are only learning. In particular, the agents density propagates at the Fisher-KPP speed. 
We name this type of solutions a   lottery society.  
\end{abstract}


%

\section{Introduction and main result}
The Lucas-Moll system  was introduced to study growth in 
macroeconomics~\cite{Lucas-Moll} through diffusion of knowledge. It describes the evolution of agents that divide their time
between production and acquiring knowledge from other agents. Their search strategy for knowledge 
is based
on the solution to an optimization problem that maximizes their overall discounted income stream. 
The original formulation in~\cite{ABLLM,Lucas-Moll} involves the agents cumulative
distribution function $F(t,x)$ and the value function $V(t,x)$. Here, $x\in\Rm$ is the 
logarithm of  the knowledge (productivity) of the individual agents. 

As observed in~\cite{Porretta-Rossi}, 
it is convenient to reformulate the Lucas-Moll system in terms of the
function 
\be
w(t,x)= {(\rho-\kappa)}e^{-x}V_x(t,x).
\ee
Here, the constants~$\kappa>0$ and $\rho>0$ 
are the diffusion coefficient (internal innovation rate in economics) and  the discount rate, respectively.  
With this definition,  
the formulation of~\cite{ABLLM,Lucas-Moll} 
is equivalent to a coupled system of forward and backward
nonlocal semilinear parabolic equations
\be\label{eqn:intro1} 
\bal
&\pdr{F(t,x)}{t}=\kappa \pdrr{F(t,x)}{x}+F(t,x)\int_{-\infty}^x\alpha(s_m(I(t,y))(-F_y)(t,y)dy,
\\
&\pdr{w(t,x)}{t}+\kappa\pdrr{w(t,x)}{x}+2\kappa w_x(t,x)+(\rho-\kappa)\big[1-s_m(I(t,x))-w(t,x)\big]\\
&~~~~~~~~~~~-\alpha(s_m(I(t,x)))w(t,x)F(t,x)=0, \\
\enbal
\ee
posed for $0\le t\le T$ and $x\in\Rm$. As~$w(t,x)$ is a weighted derivative of the value function,
we  refer to~(\ref{eqn:intro1}) as ``the derivative formulation" of the Lucas-Moll system. 

The function $I(t,x)$ is defined by  
\be\label{24sep502}
\bal
&I(t,x)=\farc{1}{\rho-\kappa}e^{-x}\int_x^\infty e^yw(t,y)F(t,y)dy,
\enbal
\ee
and we will refer to it as the learning pay-off functional, and to~$\alpha(s):[0,1]\rightarrow \Rm^+$ 
as the search for knowledge function.  
Economic considerations require $\alpha(s)$ 
to be increasing in $s$ and concave.
In addition, we will assume that 
\be\label{24sep302}
\hbox{$\alpha(0)=0$ and 
$\alpha'(0)=+\infty$.}
\ee

The function $s_m(I)$ in~\eqref{eqn:intro1} is the optimal proportion of time an individual 
allocates to searching for greater knowledge, 
instead of generating income. Intuitively,  $s_m(I)$ is the incentive to learn based on the expectation that acquiring additional knowledge will improve an individual's future income stream.
We will explain the derivation of its formula in Section~\ref{sec:origins}. Now, we just set it as follows: 
\be\label{eqn:sm}
s_m(I)=1,~~\hbox{ for }~~I>I_c:=\farc{1}{\alpha'(1)}, \hbox{ and } \alpha'(s_m(I))=\farc{1}{I},~~\hbox{for}~~I\le I_c=\farc{1}{\alpha'(1)}.
\ee
Note that $s_m(I)$ is defined uniquely by (\ref{eqn:sm}) if $\alpha(s)$ is strictly concave. 
In particular, individuals with knowledge $x\in\Rm$ such that $I(t,x)\ge I_c$ allocate all of their
time to learning, and those with~$I(t,x)<I_c$ spend a non-trivial fraction of their time
producing. 

In a sense, the function $w(t,x)$ 
is the propensity of an individual agent 
to learn. Typically, it is increasing  and takes values
\be
0\le w(t,x)\le 1.
\ee 
On the other hand, as $F(t,x)$ is a cumulative distribution function, it is decreasing and satisfies
\be
0\le F(t,x)\le 1.
\ee
We discuss in detail  the derivation of the Lucas-Moll system and its economics
interpretation
in Section~\ref{sec:origins}.  The connection between its original form in~\cite{Lucas-Moll} and
the derivative formulation (\ref{eqn:intro1}) 
is explained in Appendix~\ref{sec:append-deriv} below.

The system (\ref{eqn:intro1}) is supplemented 
with the boundary and initial conditions for~$F(t,x)$: 
\be\label{24sep306}
F(t,-\infty)=1,~~F(t,+\infty)=0,~~F(0,x)=F_0(x),
\ee
and the boundary and terminal conditions for the function $w(t,x)$:
\be\label{24sep308}
w(t,-\infty)=0,~~w(t,+\infty)=1,~~w(T,x)=w_T(x).
\ee

\subsection{Existence of traveling waves} 

Our approach to the study of the Lucas-Moll system stems from the fact 
that~\eqref{eqn:intro1} is a non-local generalization of the classical 
Fisher-KKP reaction-diffusion equation.  In fact,  in the special case~$\alpha(s)=\alpha_1=\hbox{const}$, 
the first equation in~\eqref{eqn:intro1} reduces exactly to the Fisher-KPP equation
\begin{equation}\label{feb2106}
\pdr{F}{t}=\kappa \pdrr{F }{x}+\alpha_1 F (1-F ).
\end{equation}
This equation admits traveling 
wave solutions $F(t,x)=F_*(x-ct)$ for all 
\be\label{24aug604}
c\ge c_*^{kpp}=2\sqrt{\kappa\alpha_1}.
\ee 
There is vast literature on the long time behavior of the solutions to the Fisher-KPP equations
and their convergence to traveling waves in the long time limit: see, for example,~\cite{AHR-conv}
and references therein. 
The front-like nature and linear in time spreading 
of the solutions to the Fisher-KPP equation and related models 
have been used in macroeconomics 
to model the growth of an economy since at least the 1980's~\cite{JR,KL,Luttmer1,Luttmer2,Staley}, with 
more recent
contributions, among others, in ~\cite{ABL, BPT,Ben-Bru-Hag,Ben-Per-Ton,Bue-Luc,KLZ,Lucas-Moll,Perla-Tonetti,Per-Ton-Waugh},
including an interesting recent direct comparison to data in~\cite{Staley-24}. The coupled
Lucas-Moll system (\ref{eqn:intro1}) (albeit written in the form (\ref{Psi2})-(\ref{V2}) below) 
itself was first proposed in~\cite{Lucas-Moll} without any diffusion and, simultaneously, a very much related
dynamics was introduced in~\cite{Perla-Tonetti}.  The model with diffusion was, to the best of our knowledge, first
introduced in~\cite{ABLLM}. 

The structure of the full coupled problem~\eqref{eqn:intro1}  
for a general search function~$\alpha(s)$ still inherits some
Fisher-KPP features that make one expect a similar result on the existence and long time stability of traveling wave 
solutions to~\eqref{eqn:intro1}  of the form 
\begin{equation}\label{eqn:trav}
F(t,x)=F_*(x-ct),~~w(t,x)=w_*(x-ct),
\end{equation}
known in economics as balanced growth paths. They are monotonic in $x$ and satisfy the boundary
conditions
\be
F_*(-\infty)=1,~F_*(+\infty)=0,~w_*(-\infty)=0,~w_*(+\ifnty)=1.
\ee 
At present, however, we do not have  
mathematical methods 
that allow to capitalize easily on this expectation. 
The first mathematical results on the existence of traveling waves for~\eqref{eqn:intro1} were obtained
in~\cite{BLW1,BLW2} for the case $\kappa=0$ when diffusion is absent. In that situation, the problem was analyzed from the
nonlinear Boltzmann equation point of view. 
In the presence of diffusion and for sufficiently large $\rho$, existence of the traveling waves  was proved in~\cite{PRV} for 
$\alpha(s)=\alpha_1\sqrt{s}$,
and in \cite{Porretta-Rossi} for a more general class of
$\alpha(s)$ and under sharper assumptions on the basic parameters. Both of these papers 
used a combination of 
reaction-diffusion and Hamilton-Jacobi-Bellman  equations techniques.  
\begin{thm}\label{2.1}(Porretta and Rossi~\cite{Porretta-Rossi})\\
(1) If~$\rho>\kappa$ and~$\alpha(1)>\kappa$, 
then there exists $c\in(0,2\sqrt{\kappa\alpha(1)})$ so that 
the system ~\eqref{eqn:intro1} has 
a solution of the form (\ref{eqn:trav}), such that~$F_*(x)$ is monotonically 
decreasing and~$w_*(x)$ is monotonically increasing.
In addition, there exists 
$\etal\in\Rm$ so that  $s_m(I(x))=1$ for all $x<\etal$ and
$s_m(I(x))<1$ for all $x>\etal$.  
\\
(2) Any traveling wave solution   to ~\eqref{eqn:intro1} satisfies
$2\kappa<c<\kappa+\alpha(1)$.\\
(3) There exists a traveling wave for all $c\in[2\sqrt{\kappa\alpha(1)},\kappa+\alpha(1))$.
\end{thm}
Unlike in the Fisher-KPP case, the precise range of possible traveling wave speeds 
for the Lucas-Moll system is not known but 
we expect 
that there is an interval of speeds~$[c_{min},c_{max})$ so that~\eqref{eqn:intro1} has traveling wave solutions for 
all~$c\in[c_{min},c_{max})$. That question is still open: note the gap between the speeds
in parts (1) and (3) of Theorem~\ref{2.1}. The finite range of possible speeds makes the situation very different
from the Fisher-KPP situation.

\subsection{The lottery society} 
 
The assumption 
\be\label{24sep514}
\rho>\kappa
\ee
in Theorem~\ref{2.1} that the discount rate is sufficiently large is consistent 
with economic intuition that the discount rate needs to compensate for the internal
innovation rate. We will make this assumption
throughout the paper. 

The second assumption in Theorem~\ref{2.1} that $\alpha(1)>\kappa$ means that learning is 
sufficiently efficient compared to the internal innovation. This leads to an 
equilibrium between the external learning
and internal innovation that is reflected in the existence of the balanced growth paths.
In this paper, we consider the complementary
regime~$0<\alpha(1)<\kappa$ when learning can be much less efficient than the internal innovation. In this regime,
it has been shown in~\cite{PRV,Porretta-Rossi} that balanced growth paths do not exist. Here, we are interested not in the balanced
growth paths but in the
long time behavior of the solutions to the initial-terminal value problem for the Lucas-Moll system. 
As we will see, somewhat paradoxically, even though learning is inefficient, 
this leads to a regime where agents 
spend an overwhelming fraction of their time just learning and not producing.   

To describe this situation, consider two characteristic transition points of the solutions to~\eqref{eqn:intro1}. The first point is the learning
front location $\eta_\ell(t)$, defined  via 
\be\label{24sep423}
I(t,\eta_\ell(t)) =\farc{1}{\alpha'(1)},
\ee
so that 
\be\label{24sep422}
\bal
&s_m(I(t,x))=1, 
~~\hbox{ for $x<\eta_\ell(t)$},\\
&s_m(I(t,x))<1, 
~~\hbox{ for $x>\eta_\ell(t)$.}
\enbal
\ee
That is, the agents with the productivity $x<\etal(t)$ spend all of their time learning and
do not produce. 
The second transition point is the median front $\eta_m(t)$ of the agents determined by
\be\label{eqn:median}
F(t,\eta_m(t))=\farc{1}{2}.
\ee
This is, roughly, the typical knowledge of the agents. 
In a balanced growth path, these two points are separated by a distance that is constant in time,
so that a non-trivial fraction of the agents is involved in producing.  
In contrast, it turns out that in the regime $0<\alpha(1)<\kappa$, the distance between 
the learning front location $\eta_\ell(t)$ and the median front~$\eta_m(t)$ grows in time, for general initial
and terminal conditions.  We named such solutions {\it the lottery society}, to reflect the fact that 
the agents, instead of producing, essentially gamble on a very low probability event of finding a very advanced agent that would allow them to make
a huge jump in productivity.

In other words, in a lottery society the median front is located far behind the learning front:
\be\label{eqn:fronts}
\eta_m(t)\ll \eta_\ell(t), \hbox{ and } \etal(t)-\eta_m(t)\to+\infty,~\hbox{as $t\to+\infty$.}
\ee  
In particular, the optimal proportion of time an agent is searching for new knowledge satisfies 
\be\label{eqn:search}
s_m(t,x)= 1 \hbox{ unless $F(t,x)\ll 1$.}
\ee
That is, in a lottery regime, the overwhelming majority of agents are focused solely on acquiring knowledge 
and do not engage in production, even though learning itself is rather inefficient. 
Moreover, the small fraction of agents who do produce tends 
to diminish over time as the learning and median fronts diverge~\eqref{eqn:fronts}. 

The main objective of this paper is to demonstrate that the lottery society is the long time behavior of the 
general solutions to the initial-terminal value problem for the Lucas-Moll system when $0<\alpha(1)<\kappa$. 
\begin{thm}\label{conj-feb2202_intro}
Consider solutions to~\eqref{eqn:intro1} with  $\alpha(s)=\alpha_1\sqrt{s}$, $\rho>\kappa$ and $\alpha_1< \kappa$. 
Suppose that the initial condition $F_0(x)$ is   decreasing, 
satisfies $F_0(x)=1$ for all $x\le -L_0$, and $F_0(x)=0$ for all $x\ge L_0$, and that 
 the terminal condition $w_T(x)$ is increasing and satisfies $w_T(-\infty)=0$, and~$w_T(+\infty)=1$.
Then, we have the following asymptotics for the median and learning fronts:
\be\label{as:intro}
\eta_m(t)=2\sqrt{\kappa\alpha_1}t+o(t), \hbox{ and } \etal(t)=(\kappa+\alpha_1)t+o(t).
\ee
\end{thm}
As a consequence of this result, we observe that under the assumptions of Theorem~\ref{conj-feb2202_intro}, 
the learning and agent fronts diverge, as described in~(\ref{eqn:fronts}).
Additionally, the overwhelming majority of agents are focused solely on searching and not on production, as 
indicated in~(\ref{eqn:search}).
The specific form of the assumption $\alpha(s)=\alpha_1\sqrt{s}$ is made
primarily for simplifying certain computations, 
and we anticipate that the result will hold in a broader context. On an informal level, a consequence 
of~(\ref{eqn:search}), which follows immediately from (\ref{as:intro}), is that $F(t,x)$ approximately
satisfies the Fisher-KPP equation~(\ref{feb2106}). Thus, one may expect that, in addition to the front 
asymptotics (\ref{as:intro}), the profile of $F(t,x)$ converges to that of a Fisher-KPP traveling wave as $t\to+\infty$. 
However, the proof of that convergence would likely
require a more refined asymptotics for the front location than~(\ref{as:intro}),
in the spirit of the convergence proof in~\cite{NRR1} and we do not address this question here. 

Theorem~\ref{conj-feb2202_intro} seems to be the first result on the long time behavior  of the solutions
to the Lucas-Moll system in any setting. Their behavior in the general case, without the assumption $\alpha(1)<\kappa$
remains an open question, with some results in that direction obtained in the forthcoming paper~\cite{PorRosRyz}. 
 
As we have mentioned, the conclusion that ``nearly everyone" is only learning and not
producing anything in the regime
when learning is inefficient ($\alpha_1$ is small) may seem counterintuitive at first sight.
One could expect, on the contrary, that if learning is inefficient then agents might prioritize production over learning. 
The quantitative reasons for that behavior are described in Section~\ref{sec:lottery-intro}. 
The main observation is   that 
the solution
$F(t,x)$ to the Lucas-Moll system has a relatively fat tail when $\alpha(1)<\kappa$. 
This boosts the attractiveness of learning
for the individual agents. In addition, while the chances of encountering a very advanced agent 
are small for a typical agent, the huge potential benefit outweighs the small chances.

{\bf Organization of the paper.} The paper is organized as follows. 
In Section~\ref{sec:origins} we discuss a general class of learning models
and the derivation of the Lucas-Moll system.  Section~\ref{sec:basic}
describes some basic properties of the solutions to the Lucas-Moll system.
Theorem~\ref{conj-feb2202_intro} is proved in Section~\ref{sec:lottery}.
The two appendices contain some auxiliary results.

{\bf Acknowledgement.} LR was supported by an
NSF grant DMS-2205497 and by an~ONR grant N00014-22-1-2174.  AN was supported by an
NSF grant DMS-2407046. HB and JMR have received funding from the French ANR Project ANR-23-CE40-0023-01 ReaCh. 
JMR acknowledges an invitation by the Stanford Mathematics department in winter 2024, 
that was an important step in the completion of this work.
We are grateful to Alessio Porretta and Luca Rossi for 
illuminating discussions of the Lucas-Moll system, 
and to Benjamin Moll and Christopher Tonetti for generously patient explanations of the
economics background and related work. 

\section{Origins of the problem and previous results}\label{sec:origins}

In this section we recall the origins of the Lucas-Moll system 
and discuss its connections to other systems of interacting particles. 
Generally, diffusion of knowledge type dynamics in macroeconomics attempts to provide a basic  mechanism
that leads to the growth of an economy over long time periods. 
These models  are vastly simplified but  
likely capture some of the basic features of economic growth, and also  provide 
a rich class of interesting new mathematic models.  

Before we describe the mathematical models, let us mention that to the best of our understanding, 
the Pickwickians were the first to think of the
importance of the diffusion of knowledge  and learning itself and of its connection
to production~\cite{Dickens}:
``That while this Association is deeply sensible of the advantages which must accrue to the cause of science, from the production to which they have just adverted -- 
no less than from the unwearied researches of Samuel Pickwick, Esq., G.C.M.P.C., in Hornsey, Highgate, Brixton, and Camberwell -- they cannot but entertain a lively sense of the inestimable benefits which must inevitably result from carrying the speculations of that learned man into a wider field, from extending his travels, and, consequently, enlarging his sphere of observation, to the advancement of knowledge, and the diffusion of learning." 

This section is organized as follows. In Section~\ref{sec:agents} we discuss a general class of agent based learning models. These are interacting particle
systems where individuals can jump to the locations of other agents in the ``knowledge space", as explained in Section~\ref{sec:agents-basic}.  Its formal mean-field macroscopic dynamics is discussed in
Section~\ref{sec:macro}. Diffusion (internal innovation) is added to that system in Section~\ref{sec:innovate}. Section~\ref{sec:fkpp} discusses
the connection between the resulting non-local equation and the Fisher-KPP equation, as well as an interpretation of the N-BBM process, a well-studied
stochastic system, as an example of a learning dynamics.  Section~\ref{sec:strategy} discusses some possible prescriptions of learning strategies. Until that 
point, the learning dynamics does not involve the value function or optimization in any way. Learning models based on optimization problems for the
learning strategy are discussed 
in Section~\ref{sec:value}. The Lucas-Moll model is derived as a Nash equilibrium for such learning dynamics with the search strategy
resulting from an optimization problem in Sections~\ref{sec:value-ind}, \ref{sec:value-best} and~\ref{sec:value-nash}. Finally, in Section~\ref{sec:value-fronts}
we discuss an approximation to the learning front location that is intrinsic to the particles, in the sense that it can be computed purely based on the current
particle locations, without the computation of the value function.

\subsection{The dynamics of agents}\label{sec:agents}

We first describe the basic underlying dynamics in a general class of the diffusion of learning models,
and will later specify some choices that lead to the Lucas-Moll system.  
They are similar in spirit to the Brunet-Derrida type of models 
of branching processes introduced in~\cite{Bru-Der97,Bru-Der-Mu-Mu06,Bru-Der-Mu-Mu07}.  
Consider a population of $N \gg 1$ interacting particles, which we will refer to as agents (e.g., firms). 
Each agent $A_k$ is characterized by a level of knowledge or productivity $z_k(t) \geq 0$ that varies over time. 
While the notion of knowledge as one-dimensional is simplistic -- since a company’s expertise could be highly 
multidimensional -- we will not address this complexity at present. The informal definition of productivity is that 
the value of the goods produced by an agent with current productivity $z(t) \geq 0$ during 
a time interval $[t, t + \Delta t]$ is $z(t) \Delta t$, assuming the agent focuses solely on 
production during that interval.

If an agent's goal is to maximize the value of 
its total production over time, with a proper
discounting of the future production, then it may be of a benefit not just to
solely produce while remaining at a fixed knowledge value but also to spend a fraction of time  increasing its productivity via
learning even if at the cost of temporarily
reducing the current production volume.  We will assume that each agent spends the fraction $s_k(t)\in[0,1]$ of an infinitesimal time interval 
$[t,t+dt]$ on learning and the remaining fraction $(1-s_k(t))$ on producing. Thus, its total production during this time interval
is 
\be\label{24aug602}
dP=(1-s_k(t))z_k(t)dt.
\ee
We will refer to $s_k(t)$ as the learning strategy of the agent $A_k$. 
The two basic modeling issues in the above dynamics are \\
(1)~how to model 
the acquisition of knowledge, \\
(2) how to choose the learning strategy $s_k(t)$ of each agent
so as to balance the future benefit of the acquired
knowledge and the current cost of its acquisition.

\subsubsection{Acquisition of knowledge}\label{sec:agents-basic}

Agents in the economy learn by searching for other agents with a higher production knowledge and then acquiring their 
knowledge if they meet them.  We will model this process as follows. The system consists of a large number $N\gg 1$ of 
individual agents that interact as follows. The  agents~$A_k$,~$k=1,\dots,N$ have independent exponential clocks that run at the rates
$R_k(t)$, so that the time $\tau_k$ before the clock rings has the probability law
\be
\Pm(\tau_k>t)=\exp\Big(-\int_0^t R_k(t')dt'\Big).
\ee
When a clock $\tau_k$ rings, the agent $A_k$ chooses another agent $B$, uniformly among all the other agents.  
If the productivity~$z_{A_k}(t^-)$  of  the agent $A_k$ before the clock rings is lower than the productivity~$z_B(t^-)$ of the agent~$B$, then 
the agent $A_k$ updates its production knowledge to that of the agent $B$, so that~$z_{A_k}(t^+)=z_B(t^-)$. 
If, on the other hand,~$z_{A_k}(t^-)\ge z_B(t^-)$, then
the productivity of~$A_k$ does not change, and~$z_{A_k}(t^+)=z_{A_k}(t^-)$.   In either case,
the exponential clock of the agent $A_k$ is reset at that time. Let us emphasize that the searching agent does not control which other agent 
it will find in the search. 

We will define the clock rates $R_k(t)$ in terms of the aforementioned search strategies $s_k(t)$ that appear in (\ref{24aug602}), via
\be\label{24aug502}
R_k(t)=\alpha(s_k(t)).
\ee 
Here,~$\alpha(s)$ is a given increasing function. We will assume that
it satisfies assumptions~(\ref{24sep302}). The assumption~$\alpha(0)=0$ in~(\ref{24sep302}) means that
if no time is allocated to searching then the agent can not acquire a higher knowledge. 
This assumption can sometimes be relaxed. The law of
diminishing returns implies that $\alpha(s)$ should be concave. 
The second assumption in ~(\ref{24sep302}) that
$\alpha'(0)=+\infty$ encourages the ``very elite" agents to still spend a small but
positive fraction
of time searching. 
A typical example to keep in mind is $\alpha(s)=\alpha_1\sqrt{s}$ with some $\alpha_1>0$.

An important feature in this model coming from the economic intuition 
is that the strategy $s_k(t)$ should depend
not only on the current knowledge $z_k(t)$ of a given agent
but also on the positions of the other agents. For example, one possible learning
strategy is
\be\label{24apr706}
s_k(t)=\farc{1}{N}\#\{z_m(t)\ge z_k(t),~1\le m\le N\},
\ee
which says that fraction of time an agent searches is proportional to its rank among all agents.
A smoother version of that choice is
\be\label{24apr708}
s_k(t)=\farc{1}{N}\sum_{k=1}^N\zeta(z_m(t)-z_k(t)),
\ee
with a smooth increasing function $\zeta(z)$ such that $\zeta(z)=0$ for all $z\le 0$ and $\zeta(z)=1$ for $z\ge 1$. 

Another interesting for us case is 
\be\label{24apr702}
s_k(t)=\min\Big[1,\farc{1}{N}\sum_{m:z_m(t)\ge z_k(t)}\Big(\frac{z_m(t)}{z_k(t)}-1\Big)\Big],
\ee
or, more generally, 
\be\label{24apr702bis}
s_k(t)=\min\Big[1,\farc{1}{N}\sum_{m:z_m(t)\ge z_k(t)}\Big(\farc{\vphi(z_m(t))}{\vphi(z_k(t))}-1\Big)\Big],
\ee
with a given positive increasing function $\vphi(z)$. With all of the above strategies, 
the agents that are ``lagging behind"
search more than the ``advanced agents".

A crucial modeling assumption here is that the choice of the ``knowledge giver" agent $B$ does not depend on the 
initial gap between the productivities of $A$ and $B$. In particular, 
the ``instantaneous"
jump in knowledge can be arbitrarily large if $z_B(t^-)\gg z_A(t^-)$. 
This can be restricted in
many ways but we will not focus on this issue. Let us just comment that,  
even if we relax this assumption, in order to keep the agents spreading at a positive speed, 
we will need to allow some very large jumps of productivity,
with a reasonably large probability. In other words, the tail distribution of the
jumps needs to be sufficiently heavy. 
Otherwise, if the interactions are local, then 
the system will only exhibit
diffusive spreading like $O(\sqrt{t})$, 
while in the real world economies grow at a constant rate in time.

\subsubsection{The macroscopic dynamics}\label{sec:macro}

Let us emphasize that even though the above probabilistic descriptions  of the models  
are discussed for a finite~$N\gg 1$, the Lucas-Moll system 
is obtained in  a ``mean-field" macroscopic description, in the limit $N\to+\infty$.
The passage to the mean-field limit is, in itself, an interesting open problem. 
 
Let us now describe the macroscopic dynamics that arise formally in the limit as $N\to+\infty$. 
From this point forward, we will work with the logarithmic variable $x=\log z$. In the following discussion, we will often refer to $x$
as knowledge or productivity, even though it represents log-productivity rather than actual productivity.
We let $\psi(t,x)$ be the probability density of the agents at a time~$t\ge 0$ 
obtained in the limit $N\to+\ifnty$.  Then, 
the probability that the agent~$B$, who~$A$ encounters, has log-knowledge 
in the interval~$(x',x'+\Delta x')$
is~$\psi(t,x')\Delta x'$.    
We will assume that the agents use
a common search strategy~$s_c(t,x)$ that depends only on the log-productivity~$x$ and time~$t\ge 0$, 
but not on an individual agent.  Once $s_c(t,x)$ 
 is specified, the evolution of the agent density is governed by an integral equation
\begin{equation}\label{24apr402}
\begin{aligned}
\pdr{\psi(t,x)}{t}&= 
\int_{-\infty}^x\alpha(s_c(t,y))\psi(t,x)\psi(t,y)dy
- \int_x^\infty \alpha(s_c(t,x))\psi(t,x)\psi(t,y)dy.
\end{aligned}
\end{equation}
The first term in the right
side  is due
to the agents  at knowledge $y\le x$ who execute a successful search and end up
at knowledge $x$, and the second one comes from the agents with knowledge~$x$ 
who do a successful search and move up to a level $y\ge x$.

%

\subsubsection{Innovation and diffusion of knowledge} \label{sec:innovate}

Another important component of the dymamics  is that,
independent of the learning from other agents, each agent also attempts to
innovate internally, a process that has nothing to do with the distribution
of the other agents. Innovation can lead both to an improvement and a regression, and is
often modeled by a Brownian motion. Once again, for simplicity, we will disregard the
cost of innovation, though some of it may be attributed to the "regression due to innovation"
inherent in the Brownian motion. From the economics point of view, it is reasonable to assume that the "innovation diffusion"  takes place 
in the logarithmic variable $x=\log z$ and not in the productivity variable $z$ itself.  

Note that this diffusion via innovation is very different from
the Pickwickian diffusion of knowledge~\cite{Dickens}: 
here, diffusion is an internal process for each given
agent. On the other hand, in the
Pickwickian model, diffusion of one agent is crucial for other agents:  
there, an agent performs a Brownian motion (note that Brownian motion is an excellent model for
the travels of Mr.~Pickwick) and transmits knowledge to  
agents it encounters on the way.  The Pickwickian model is also an interesting interacting particle system
but is very different from the diffusion of knowledge models we consider here.

The presence of innovation modifies  
the macroscopic evolution  equation (\ref{24apr402})  for the common
agent density~$\psi(t,x)$ to
\begin{equation}\label{may2106}
\begin{aligned}
\pdr{\psi(t,x)}{t}&=\kappa \pdrr{\psi(t,x)}{x}+
\int_{-\infty}^x\alpha(s_c(t,y))\psi(t,x)\psi(t,y)dy
- \int_x^\infty \alpha(s_c(t,x))\psi(t,x)\psi(t,y)dy,
\end{aligned}
\end{equation}
with a prescribed initial condition $\psi(0,x)=\psi_0(x)$. 

 \subsubsection{Evolution of the cumulative distribution function} \label{sec:fkpp}

Equation (\ref{may2106}) has the form of a nonlinear Boltzmann equation from the 
kinetic theory with an extra diffusion term. 
A very different point of view is based on the evolution of  the cumulative distribution function
\begin{equation}\label{feb2102}
F(t,x)=\int_{x}^\infty\psi(t,z)dz.
\end{equation}
It satisfies a non-local  equation
\begin{equation}\label{feb2104}
\pdr{F(t,x)}{t}=\kappa \pdrr{F(t,x)}{x}+F(t,x)\int_{-\infty}^x \alpha(s_c(t,z))(-F_z(t,z))dz,
\end{equation}
with the boundary conditions $F(t,-\infty)=1$, $F(t,+\infty)=0$.

As we have mentioned, in the special case $\alpha(s)=\alpha_1=\hbox{const}$, this equation
reduces exactly to the classical Fisher-KPP equation
\begin{equation}\label{feb2106bis}
\pdr{F}{t}=\kappa \pdrr{F }{x}+\alpha_1 F (1-F ).
\end{equation}
At the microscopic level, this is the case when the exponential clocks of all agents are identical and do not
depend on their positions or time: $R_k(t)\equiv \alpha_1$. In particular, the only interaction between the agents
is through them jumping on top of each other but the search strategy of an agent does not
depend on the location of the other agents.

When the search strategy of each agent does depend on the locations of the
other agents, the learning process described above has a similarity with a branching process known as $N$-BBM. 
This is a version of 
a branching Brownian motion on the real line, introduced by Brunet, Derrida, Mueller and Munier in~\cite{Bru-Der-Mu-Mu06,Bru-Der-Mu-Mu07}
and later studied in~\cite{Berard-Mai,Ber-Zhao,Ber-Bru-Pen,Ber-Bru-Nol-Pen,Ber-Bru-Nol-Pen-Bees,Dur-Rem}. 
That process also consists of $N$ particles that may diffuse and have exponential clocks attached to them. When the clock rings,
the corresponding particle branches into two and, simultaneously, the particle with the smallest value of a prescribed scoring function 
is removed. The overall number $N\gg 1$ of particles stays constant in the dynamics. In one dimension
a natural scoring function is simply $s(x)=x$. In higher dimensions, 
the cases~$s(x)=\|x\|$ and $s(x)=1/\|x\|$ have been considered. In the latter case, the
system is known as Brownian bees. 

In the present  context, the $N$-BBM corresponds 
to a version of learning dynamics 
where only the  agent with the smallest knowledge
is allowed to learn. More generally, one 
can think of the learning process we described above as the following version of $N$-BBM.
Suppose that in the N-BBM process, when a particle branches, instead of removing from the system the particle with the smallest value of the scoring
function, we do the following. 
First, another particle is chosen uniformly 
among all particles. Then, we compare the values of the scoring function of the newly born particle and the chosen particle.  
The particle with the smaller scoring function between these two particles 
is removed from the system, so that the total number of particles stays fixed. 
Unlike in N-BBM, the value of the scoring function may depend not only on the
location of a given particle but also on the locations of all other particles.

One of the most interesting properties of $N$-BBM 
was proved  by N. Berestycki and L. Zhao~\cite{Ber-Zhao}. They proved that in 
dimensions~$d>1$ the~$N$-BBM particles with the scoring
function~$\|x\|$ move, as a clump, with a deterministic speed but in a random direction chosen
uniformly on the unit sphere. It would be interesting to show a similar phenomenon 
in the multi-dimensional versions
of the learning process considered in the present paper. In particular, when the productivity
of the firms is considered as a multi-dimensional object, this may lead to 
a choice of a random direction in which the vast majority of the firms would be moving.


\subsubsection{Possible prescriptions of a strategy}\label{sec:strategy} 

Let us now discuss what would constitute a sensible choice of the
search strategy~$s_c(t,x)$.  The very unproductive agents would be learning almost all the time, or even
all the time. On the other hand,
the very advanced agents would spend very little time learning because for them
the chances of meeting a more knowledgeable agent are very small. Internal innovation may 
give them a higher chance of a technological advance than an outside search. 
Hence, the search strategy~$s_c(t,x)$ should have the limits
\be\label{24apr504}
s_c(t,-\infty)=1,~~s_c(t,+\infty)=0.
\ee
To capture the above intuition, one natural
choice is to postulate that the common search strategy is given by 
\be\label{may2902bis}
s_c(t,x)=F(t,x).
\ee
This corresponds to the search strategy
we have seen in (\ref{24apr706}) for a finite $N\gg 1$: the search strategy of an agent
depends solely on its rank.  
In that case,
the non-local Fisher-KPP equation~(\ref{feb2104}) becomes local:
\begin{equation}\label{may2202}
\pdr{F}{t}=\kappa \pdrr{F}{x}+F\big(Q(1)-Q(F)),
\ee
with the nonlinearity 
\be
Q(u)=\int_0^u\alpha(u)du.
\end{equation}
Equation (\ref{may2202}) is a reaction-diffusion equation
of the Fisher-KPP type and its solutions converge to a traveling wave as $t\to+\infty$. 
However, (\ref{may2902bis})  
is the only known choice of the common search strategy $s_c(t,x)$
when such result on spreading is readily available. Even for the slightly non-local search strategies
\be\label{24apr710}
s_c(t,x)=\int_\Rm F(t,x-y)\zeta(y)dy,
\ee
with a nice rapidly decaying kernel $\zeta(x)\ge 0$, with $\|\zeta\|_{L^1}=1$, 
the long time convergence 
of the solutions to (\ref{feb2104}) to a traveling wave is not known. This choice
of strategy corresponds to (\ref{24apr708}) for a finite $N\gg 1$.

\subsection{The value function and its evolution} \label{sec:value}

So far, we have described the dynamics of the density of agents assuming that the
common search strategy $s_c(t,x)$ is given. In macroeconomics models, the 
search strategies of interest come from optimization problems that we now describe.


\subsubsection{The value function for a given individual strategy}\label{sec:value-ind} 

Let us first suppose that ``nearly all" agents use a common search 
strategy $s_c(t,x)$ and their cumulative distribution function 
satisfies
the non-local Fisher-KPP equation (\ref{feb2104}), but a ``black market" 
individual agent   follows
an individual strategy $s_i(t,x)$ that may differ from $s_c(t,x)$. Then, 
the cumulative distribution function $G(\tau,y;t,x)$ for the individual agent satisfies
a linearized Fisher-KPP type equation
\begin{equation}\label{feb2118}
\pdr{G(\tau,y;t,x)}{\tau}=\kappa \pdrr{G(\tau,y;t,x)}{y}+F(t,y;s_c)\int_{-\infty}^y \alpha(s_i(\tau,z))(-G_z(\tau,z;t,x))dz,
\end{equation}
for $\tau>t$, with the initial condition at the time $\tau=t$
\begin{equation}\label{feb2822}
\begin{aligned}
G(\tau=t,y;t,x;s_c,s_i)&=1\hbox{ for $y<x$}, \\
G(\tau=t,y;t,x;s_c,s_i)&=0\hbox{ for~$y>x$. }
\end{aligned}
\end{equation}
If the two search strategies are identical: 
$s_i(t,x)=s_c(t,x)$ then the distributions of the common agents and the 
black market agent will, of course, coincide.

Given a common  strategy $s_c(t,x)$ and an individual strategy $s_i(t,x)$, 
we can  find (i) the common density of agents~$\psi(t,x)$, from the solution to the non-local
Fisher-KPP equation  (\ref{feb2102})-(\ref{feb2104}), 
and (ii) the individual transition probability density 
$p(\tau,y;t,x)$, from the solution
to (\ref{feb2118}).
With these in hand, we define the average individual production
corresponding to these strategies as
\begin{equation}\label{feb2016}
\begin{aligned}
Z(t,x;s_c,s_i)&=\int_t^T\int e^{-\rho(\tau-t)}e^{y}[1-s_i(\tau,y)]p(\tau,y;t,x)
dy d\tau+\int e^{-\rho(T-t)}V_T(y)p(T,y;t,x)dy.
\end{aligned}
\end{equation}  
Here, the first term is the total production between the starting time $t\ge 0$
and the fixed terminal time $T>t$, and $V_T(x)$ is a prescribed terminal condition.
The parameter $\rho>0$ represents the discounting rate. 

A  computation shows that
the function $Z(t,x;s_c,s_i)$ satisfies a linear integro-differential 
evolution equation
\begin{equation}\label{feb2018}
\begin{aligned}
\rho Z(t,x)&=\pdr{Z(t,x)}{t}+\kappa \pdrr{Z(t,x)}{x}+e^x(1-s_i(t,x))+
\alpha(s_i(t,x))\int_x^\infty \psi(t,z;s_c)[Z(t,z)-Z(t,x)]dz,
\end{aligned}
\end{equation}
with the terminal condition $Z(T,x)=V_T(x)$. Its evolution depends both on the
common strategy~$s_c$, via the common density $\psi(t,z;s_c)$ and on the
individual strategy $s_i$.

 \subsubsection{The best individual strategy for a given common strategy}\label{sec:value-best}

Let us fix a given common 
strategy~$s_c(t,x)$ followed by the vast majority of agents 
and look for an optimal strategy~$s_i(t,x)$ for an individual agent.
The agent's goal is to maximize its total production~$Z(t,x)$ given by~(\ref{feb2016}). 
In other words, we define the individual 
value function corresponding to a collective strategy $s_c$ as 
\begin{equation}\label{feb2826}
\begin{aligned}
W(t,x;s_c)=\sup_{s_i\in {\cal A}_T}
Z(t,x;s_c,s_i).
\end{aligned}
\end{equation}  
Here, ${\cal A}_T$ is the set of all admissible strategies on the time interval $[0,T]$. 
This function satisfies the Hamilton-Jacobi-Bellman equation
\begin{equation}\label{may2404}
\begin{aligned}
\rho W(t,x;s_c)&=\pdr{W(t,x;s_c)}{t}+\kappa \pdrr{W(t,x;s_c)}{x}\\
&+\max_{s\in[0,1]}\Big(e^x(1-s) +
\alpha(s)\int_x^\infty \psi(t,z;s_c)[W(t,z;s_c)-W(t,x;s_c)]dz\Big),
\end{aligned}
\end{equation}
with the terminal condition $W(T,x;s_c)=V_T(x)$.  

Given a common strategy $s_c\in {\cal A}_T$, we may now define the best individual
strategy as the maximizer in (\ref{may2404}):
\begin{equation}\label{may2406}
s_i^*(t,x;s_c)=\hbox{argmax}_{s\in[0,1]}\Big(e^x(1-s) +
\alpha(s)\int_x^\infty \psi(t,z;s_c)[W(t,z;s_c)-W(t,x;s_c)]dz\Big).
\end{equation} 
This gives us a map
\begin{equation}\label{may2402}
{\cal F}_T:~{\cal A}_T\to {\cal A}_T,~~{\cal F}[s_c](t,x)=s_i^*(t,x;s_c).
\end{equation}
The map ${\cal F}_T$ depends on the terminal condition $V_T(x)$ and on the initial condition
$\psi_0(x)$.

\subsubsection{The Nash equilibrium strategy and the Lucas-Moll system} \label{sec:value-nash}

If the optimal individual 
strategy $s_i(t,x)$ found in the above optimization procedure 
does not coincide with the common
strategy $s_c(t,x)$, then it is natural to expect that all agents will 
``go to the black market" and start shifting
to the best individual strategy~$s_i(t,x)\neq s_c(t,x)$. Thus, $s_c(t,x)$ will no longer be the strategy
used by the ``vast majority" of the agents. This, in turn, will  modify the
best individual strategy $s_i(t,x)$, and so on. One expects that this adjustment
would eventually result
in the system being in a Nash equilibrium -- this is a state where, given the common strategy $s_c(t,x)$, 
the best individual strategy is also $s_c(t,x)$. That is, a Nash equilibrium 
is a fixed point of the map~${\cal F}_T$ defined by (\ref{may2402}):
\begin{equation}\label{may2402bis}
{\cal F}[s_{Nash}]=s_{Nash} .
\end{equation}
This relation encodes the Lucas-Moll system. 

The fixed point equation (\ref{may2402bis}) is rather implicit. 
Explicitly, it takes the following form.
The cumulative probability distribution function of the
agents satisfies the evolution equation (\ref{feb2104})  
\begin{equation}\label{Psi2}
\begin{aligned}
\pdr{F}{t}-\kappa\pdrr{F}{x}&= F(t,x)
\int_{-\infty}^x\alpha(s^*(t,y))\big(-F_y(t,y)\big)dy,
\end{aligned}
\end{equation}
with the strategy $s^*(t,x)$ that comes from the 
Hamilton-Jacobi-Bellman equation (\ref{may2404}) for the value function $V(t,x)$
\begin{equation}\label{V2}
\rho V(t,x)=\pdr{V(t,x)}{t}+\kappa\pdrr{V(t,x)}{x}+\max_{s\in[0,1]}\Big[(1-s)e^x+\alpha(s)\int_x^\infty [V(t,y)-V(t,x)](-F_y(t,y))dy\Big].
\end{equation} 
The strategy  $s^*(t,x)$ that appears in (\ref{Psi2}) is the maximizer in (\ref{V2}), which couples the 
two equations. 

To connect the form (\ref{Psi2})-(\ref{V2}) of the Lucas-Moll system 
to the system (\ref{eqn:intro1}), 
let us  define 
\be
w(t,x)= {(\rho-\kappa)}e^{-x}V_x(t,x), 
\ee
and
\be
I(t,x)=\farc{1}{\rho-\kappa}e^{-x}\int_x^\infty e^yw(t,y)F(t,y)dy.
\ee
The maximizer in (\ref{V2}) can be computed explicitly, see Appendix~\ref{sec:append-deriv},  
and is given by 
\be\label{24sep420}
s^*(t,x)=s_m(I(t,x)).
\ee
This is why the function $s_m(I)$, defined in (\ref{eqn:sm}), appears in the Lucas-Moll system
(\ref{eqn:intro1}). 
As also shown in Appendix~\ref{sec:append-deriv} and originally observed in~\cite{Porretta-Rossi}, 
the function $w(t,x)$ 
satisfies a backward in time parabolic equation 
\be\label{24sep414}
\bal
\pdr{w(t,x)}{t}&+\kappa\pdrr{w(t,x)}{x}+2\kappa w_x(t,x)+(\rho-\kappa)\big[1-s_m(I(t,x))-w(t,x)\big]\\
&-\alpha(s_m(I(t,x)))w(t,x)F(t,x)=0.
\enbal
\ee
This is the second equation in (\ref{eqn:intro1}).  

A key feature of the system (\ref{Psi2})-(\ref{V2}) is that the density 
equation (\ref{Psi2}) is run forward in time: we need to prescribe an initial
condition $F(0,x)=F_0(x)$. On the other hand, the HJB equation~(\ref{V2})
is run backward in time: we prescribe the terminal value $V(T,x)=V_T(x)$. 
This is a common feature for mean-field games problems that makes the study
of their long time behavior extremely challenging. 

\subsubsection{The learning and intrinsic learning fronts} \label{sec:value-fronts}  
  
Let us briefly  discuss two critical locations for the Lucas-Moll system. The first  one is 
the learning front location defined by (\ref{24sep423}).  To the left of it, the search
strategy is 
\be
\hbox{$s^*(t,x)=s_m(I(t,x))=1$ for~$x<\etal(t)$,}
\ee
as in~(\ref{24sep422}), 
so that these agents are only searching and not
producing. The location~$\eta_\ell(t)$ plays quite an important role as this is where the definition of $s_m(I(t,x))$
changes from $s_m(I(t,x))=1$ for $x<\etal(t)$ to $s_m(I(t,x))<1$ for $x>\etal(t)$. 
As its definition involves~$I(t,x)$, it depends both on the cumulative distribution  function $F(t,x)$ and on the value 
function $V(t,x)$. An interesting alternative to $\etal(t)$ is the intrinsic front location $\el(t)$ defined purely in terms
of~$F(t,x)$ as follows. We set
\be
J(t,x)=\farc{1}{\rho-\kappa}e^{-x}\int_x^\infty e^yF(t,y)dy.
\ee
This is essentially equivalent to replacing $w(t,y)$ in the definition (\ref{24sep502})
of $I(t,x)$ by a step function that has a sharp transition from $0$ on the left
to $1$ on the right. Then, we define $\el(t)$, similarly to~(\ref{24sep423}),  by
\be\label{24sep510}
J(t,\el(t))=\farc{1}{\alpha'(1)}. 
\ee
The advantage of $\el(t)$ over $\etal(t)$ is that this object can be computed  purely in terms of the agents locations,
and does not involve the value function.  As we will see in Sections~\ref{sec:bounds-right} and~\ref{sec:tight-learning}, the learning front is,
in a sense, tight: the width over which the learning strategy drops from $s^*(t,x)=1$ to very small values, is bounded in time, and the function
$w(t,x)$ is not that far from being a step function centered at $\etal(t)$:
\be
w(t,x)\approx 1,~~x>\etal(t).
\ee
For the value function itself, this corresponds to  the approximation 
\be\label{24sep504}
V(t,x)\approx \bar V(t,x):=\farc{1}{\rho-\kappa}e^x,~~x>\etal(t).
\ee
These estimates will be used in Section~\ref{sec:learning-intrinsic} to show that the intrinsic front location gives a good approximation
to the true learning front location. Based on this observation, it is tempting to consider a purely forward in time version of the Lucas-Moll
system:
\be\label{24sep424}
\bal
&\pdr{F(t,x)}{t}=\kappa \pdrr{F(t,x)}{x}+F(t,x)\int_{-\infty}^x\alpha(s_m(J(t,y))(-F_y)(t,y)dy,
\\
&J(t,x)=\farc{1}{\rho-\kappa}e^{-x}\int_x^\infty e^yF(t,y)dy.
\enbal
\ee
This corresponds to the learning dynamics where the agents compute their learning strategy $s^*(t,x)$ 
not based on genuinely optimizing the true value function
as in (\ref{V2}) but on the assumption that the value function $V(t,x)$ can be replaced by its approximation (\ref{24sep504}) and setting
\be
s^*(t,x)=\argmax_{x\in[0,1]}\Big[(1-s)e^x+ \alpha(s)\int_x^\infty [\bar V(t,y)-\bar V(t,x)](-F_y(t,y))dy\Big].
\ee
An  advantage of  this learning model is that it allows to compute the learning strategy directly based on the positions of the other agents
and does not involve solving the backward in time Hamilton-Jacobi-Bellman equation. 
An interesting open question is to understand how closely the solutions to the approximate model~(\ref{24sep424})
resemble those of the full Lucas-Moll system (\ref{eqn:intro1}).

\section{Basic properties of the solutions to the Lucas-Moll system}\label{sec:basic}

In this section, we discuss some of the basic properties of 
the solutions to the Lucas-Moll system. We will mostly work with its
derivative formulation (\ref{eqn:intro1}):
\be\label{24sep402} 
\bal
&\pdr{F(t,x)}{t}=\kappa \pdrr{F(t,x)}{x}+F(t,x)\int_{-\infty}^x\alpha(s_m(I(t,y))(-F_y)(t,y)dy,\\
&\pdr{w(t,x)}{t}+\kappa\pdrr{w(t,x)}{x}+2\kappa w_x(t,x)+(\rho-\kappa)\big[1-s_m(I(t,x))-w(t,x)\big]\\
&~~~~~~~~~~~-\alpha(s_m(I(t,x)))w(t,x)F(t,x)=0,
\\
&I(t,x)=\farc{1}{\rho-\kappa}e^{-x}\int_x^\infty e^yw(t,y)F(t,y)dy.
\enbal
\ee
To recapitulate, the function $F(t,x)$ satisfies the boundary conditions (\ref{24sep306}):
\be\label{24sep406}
F(t,-\infty)=1,~~F(t,+\infty)=0.
\ee
while $w(t,x)$ satisfies the boundary conditions in (\ref{24sep308}): 
\be\label{24sep408}
w(t,-\infty)=0,~~w(t,+\infty)=1.
\ee
The
learning front location $\eta_\ell(t)$ is determined by  
\be\label{24sep410}
I(t,\eta_\ell(t)) =\farc{1}{\alpha'(1)},
\ee
and the function $s_m(I)$ is defined 
by (\ref{eqn:sm}).


Let us briefly discuss the initial condition for the function $F(t,x)$ and 
the terminal condition for the function~$w(t,x)$ that should be added to (\ref{Psi2})
and (\ref{24sep414}). 
We will always suppose that 
the initial condition
\be\label{24apr1208}
F(0,x)=F_0(x)
\ee
is decreasing, and satisfies the boundary conditions~
\be\label{24jan2322}
\hbox{$F_0(x)=1$ for all $x\le -L_0$ and $F_0(x)=0$ for all $x\ge L_0$,}
\ee
with some $L_0\ge 0$.

%
We also assume that the terminal condition 
\be\label{24sep506}
w(T,x)=w_T(x)
\ee
for $w(t,x)$ is increasing and satisfies
\be\label{24jan3016}
w_T(-\infty)=0,~~w_T(+\infty)=1.
\ee

Our main interest in this section is in the properties of the solutions that
do not depend on the terminal condition $w_T(x)$ in any significant way.  In Section~\ref{sec:monot}, we show that
if the initial condition~$F_0(x)$ in (\ref{24apr1208}) is decreasing and the terminal condition $w_T(x)$ in (\ref{24sep506}) is increasing, 
then~$F(t,x)$ is decreasing in $x$, $w(t,x)$ is increasing in $x$,  and the optimal strategy $s^*(t,x)$
is decreasing in $x$,  for all $0\le t\le T$, as expected from the economic intuition. The exponential decay bounds on $s^*(t,x)$ to the right of the learning
front $\etal(t)$ are proved in Section~\ref{sec:bounds-right}. In a sense, they given informal justification of the restriction that the learning
strategy has to be a step function adopted in~\cite{Perla-Tonetti}. 
Section~\ref{sec:regular} contains the proof of some basic regularity estimates
on the solutions to the Lucas-Moll system. These estimates lead to a uniform in time bound on the instantaneous speed of the
learning front in Section~\ref{sec:bound-speed}. This, in turn, leads to the proof of tightness of the front of $w(t,x)$ around
the learning front location $\etal(t)$ in Section~\ref{sec:tight-learning}. The corresponding tightness of the front of $F(t,x)$ around the median
front location $\eta_m(t)$ is proved in Section~\ref{sec:tight-particles}.  Finally, the approximation of the learning front
location $\etal(t)$ by the intrinsic front location $\el(t)$ defined in (\ref{24sep510}) is proved in Section~\ref{sec:learning-intrinsic}. 
In particular, it gives a universal lower bound on the learning front spreading speed that depends only on the diffusion coefficient
(internal innovation rate), described in Corollary~\ref{cor-24jan1302}. 
All of these ingredients will be used later in Section~\ref{sec:lottery} in the proof of Theorem~\ref{conj-feb2202_intro}.

\subsection{Monotonicity of the solutions}\label{sec:monot} 
 
We will use extensively the monotonicity of the solutions
to (\ref{24sep402}). To this end, we now show that 
any solution to (\ref{24sep402})  with a front-like initial condition $F_0(x)$ and a front-like terminal condition~$w_T(x)$
will have the front-like structure for all $0\le t\le T$. 
\begin{prop}\label{prop-24jan1202}
Let $F(t,x)$, $w(t,x)$ be a solution to (\ref{24sep402}) for $0\le t\le T$, and 
suppose that the initial condition $F(0,x)=F_0(x)$ and the 
terminal condition $w(T,x)=w_T(x)$ 
are monotonic and, respectively, satisfy the boundary conditions  (\ref{24jan2322}) and 
(\ref{24jan3016}).
Then, we have
\be\label{23dec2710}
\bal
&F_x(t,x)<0,~~\hbox{for all $0<t<T$ and all $x\in\Rm$,}\\
&w(t,x)>0,~~\hbox{for all $0<t<T$ and all $x\in\Rm$,}\\
&w_x(t,x)>0,~~\hbox{for all $0<t<T$ and all $x\in\Rm$,}
\enbal
\ee
As a consequence, we have
\be\label{24feb518}
0\le w(t,x)\le 1,~~\hbox{ for all $0\le t\le T$ and $x\in\Rm$.}
\ee
Moreover, the optimal strategy $s^*(t,x)=s_m(I(t,x))$ is decreasing in $x$:
\be\label{23dec2712}
s_x^*(t,x)\le 0,~~\hbox{for all $0<t<T$ and all $x\in\Rm$.}
\ee
\end{prop}
{\bf Proof.} It follows immediately from differentiating the first equation in (\ref{24sep402})
that if $F_0(x)$ is decreasing, then 
\be\label{24jan1504}
\hbox{$F_x(t,x)<0$ for all $t>0$ and $x\in\Rm$. }
\ee
In addition, as $0\le s^*(t,x)\le 1$,
it follows from the second equation in (\ref{24sep402}) and the positivity of the terminal condition
$w(T,x)$ that  $w(t,x)\ge 0$ for 
all $x\in\Rm$ and $0\le t\le T$. 

As $w(t,x)\ge 0$ and $F(t,x)\ge 0$, the function $I(t,x)$ is monotonically decreasing in $x$. 
Recall also that the function~$s_m(I)$ is monotonically increasing
in $I$. It follows that that $s^*(t,x)=s_m(I(t,x))$ is  decreasing in $x$.  

To see that $w(t,x)$ is increasing, we differentiate the second equation in~(\ref{24sep402}), to obtain
the following equation for $w'(t,x)=w_x(t,x)$: 
\be
\bal
&\pdr{w'}{t}+\kappa\pdrr{w'}{x}+2\kappa w'_x(t,x)-(\rho-\kappa)\big[s_m'(I)I_x+w'\big]
\\
&-\alpha'(s_m)s_m'(I(t,x))I_xw(t,x)F(t,x)-\alpha(s_m(I))Fw'-\alpha(s_m(I))wF_x=0,
\enbal
\ee
which can be written
as
\be
\bal
&\pdr{w'}{t}+\kappa\pdrr{w'}{x}+2\kappa w'_x(t,x)-A(x)w'+B(x)=0,
\enbal
\ee
with
\be
A(x)=(\rho-\kappa)+\alpha(s_m(I))F,
\ee
and
\be
B(x)=-(\rho-\kappa)s_m'(I)I_x-\alpha'(s_m)s_m'(I(t,x))I_xw(t,x)F(t,x)-\alpha(s_m(I))wF_x.
\label{eq:bx}
\ee
The first term in the right side of~\eqref{eq:bx} 
is positive because $s_m'(I)\ge 0$ and $I_x\le 0$. The second term 
is positive because
$\alpha'(s)>0$, $s_m'(I)>0$ and $I_x<0$. Finally, the last term is positive since~$F_x<0$. As the terminal condition $w'(T,x)$ is also positive,
we conclude that $w'(t,x)$ is positive for all $0\le t\le T$.~$\Box$


\subsection{Bounds to the right of the learning front}\label{sec:bounds-right}

We have the following decay estimates to the right of the learning front $\eta_\ell(t)$. 
\begin{lem}\label{lem-mar502-right}
Assume that $\alpha(s)$ is concave with $\alpha(0)=0$, and $\alpha'(0)=+\infty$,
and, in addition, that the function $\alpha^2(s)$ is convex. Let $\eta_\ell(t)$ be determined by (\ref{24sep410}):
\be\label{24jan3018bis}
I(t,\eta_\ell(t))=\farc{1}{\alpha'(1)}.
\ee 
Then, we have the following bounds:
\be\label{24jan306}
\bal
I(t,x)&\le 
\farc{1}{\alpha'(1)}e^{-(x-\eta_\ell(t))},~~\hbox{ for all $x>\eta_\ell(t)$,}
\enbal
\ee
as well as 
\begin{equation}\label{23dec2735}
\alpha(s_m(I(t,x)))\le \alpha(1)e^{-(x-\eta_\ell(t))} ,\hbox{ for all $x>\eta_\ell(t)$.}
\end{equation}
In addition, if the function $\beta(s)=\sqrt{s}\alpha'(s)$ is increasing, then we also have 
\be\label{24jan308}
s_m(I(t,x))\le e^{-2(x-\eta_\ell(t))},~~
\hbox{ for all $x>\eta_\ell(t)$.}
\ee
\end{lem} 
The assumptions on $\alpha(s)$ in the above lemma are satisfied, for example, for 
\be\label{24jan3022}
\alpha(s)=\alpha_1s^k,~~\hbox{ with $1/2\le k<1$,}
\ee
considered in Appendix~\ref{sec:alpha-powers}. 

{\bf Proof.}  Note that for $x>\eta_\ell(t)$ we have
\be\label{24apr1214}
\bal
I(t,x)&=e^{-x}\int_x^\infty e^yw(t,y)F(t,y)dy\le e^{-x}\int_{\eta_\ell(t)}^\infty e^yw(t,y)F(t,y)dy= e^{-x}e^{\eta_\ell(t)}I(t,\eta_\ell(t))\\
&=
\farc{1}{\alpha'(1)}e^{-(x-\eta_\ell(t))}.
\enbal
\ee
As we have noted above, the function $s_m(I)$ is increasing. 
Therefore, we have from (\ref{24apr1214}) that 
\be
\bal
s_m(I(t,x))\le s_m\Big(\farc{1}{\alpha'(1)}e^{-(x-\eta_\ell(t))}\Big),~~\hbox{ for all $x>\eta_\ell(t)$.} 
\enbal
\ee
Thus, to prove (\ref{23dec2735}), it suffices to show that
\be\label{23dec2739}
\alpha\Big(s_m\Big(\farc{1}{\alpha'(1)}e^{-(x-\eta_\ell(t))}\Big)\Big)\le\alpha(1)e^{-(x-\eta_\ell(t))},~~\hbox{for $x>\eta_\ell(t)$,}
\ee
or, equivalently, that for all $0<\xi<1$ we have
\be\label{23dec2740}
\alpha\Big(s_m\Big(\farc{1}{\alpha'(1)}\xi\Big)\Big)\le\alpha(1)\xi.
\ee
This, in turn, is equivalent to 
\be\label{23dec2741}
s_m\Big(\farc{1}{\alpha'(1)}\xi\Big)\le\alpha^{-1}(\alpha(1)\xi).
\ee
As $\alpha'(s)$ is decreasing, applying $\alpha'$ to both sides of (\ref{23dec2741}), we see that it is equivalent to
\be\label{23dec2742}
\alpha'\Big(s_m\Big(\farc{1}{\alpha'(1)}\xi\Big)\Big)\ge\alpha'\Big(\alpha^{-1}(\alpha(1)\xi)\Big).
\ee
Using the definition (\ref{eqn:sm}) of $s_m(I)$, this is equivalent to 
\be
\farc{\alpha'(1)}{\xi}\ge \alpha'\Big(\alpha^{-1}(\alpha(1)\xi)\Big),
\ee
or
\begin{equation}\label{23dec2743}
\farc{1}{\alpha'(1)}  \xi \le \farc{1}{\alpha'\big(\alpha^{-1}(\alpha(1)\xi)\big)}.
\end{equation} 
Let us denote the inverse $\zeta(r)=\alpha^{-1}(r)$, so that
\[
\zeta'(r)=\farc{1}{\alpha'(\alpha^{-1}(r))}.
\]
Then, the left side of (\ref{23dec2743}) is
\begin{equation}\label{23dec2748}
\farc{1}{\alpha'(1)}\xi =\farc{1}{\alpha'(\alpha^{-1}(\alpha(1)))} \xi=\zeta'(\alpha(1))\xi,
\end{equation} 
and its right side is 
\begin{equation}\label{23dec2744}
\farc{1}{\alpha'\big(\alpha^{-1}(\alpha(1)\xi)\big)}= {\zeta'(\alpha(1)\xi)}.
\end{equation} 
Therefore, (\ref{23dec2743}) simply says that
\begin{equation}\label{23dec2746}
\zeta'(\alpha(1))\xi\le  {\zeta'(\alpha(1)\xi)},
\end{equation} 
for all $\xi\in(0,1)$. Let us set $r=\alpha(1)\xi$, then (\ref{23dec2746}) takes the form
\begin{equation}\label{23dec2745}
\farc{\zeta'(\alpha(1))}{\alpha(1)}\le  \farc{\zeta'(r)}{r},\hbox{ for all $0<r\le \alpha(1)$.}
\end{equation} 
Hence, for (\ref{23dec2735}) to hold, it suffices for the function  
\begin{equation}\label{23dec2749}
\farc{\zeta'(r)}{r}=\farc{1}{r\alpha'(\alpha^{-1}(r))}
\end{equation}
to be decreasing in $r\in(0,\alpha(1))$. Equivalently, taking $r=\alpha(s)$, we want the function 
\begin{equation}\label{23dec2750}
\farc{1}{\alpha(s)\alpha'(s)}
\end{equation}
to be decreasing in $s\in(0,1)$. This, in turn, is equivalent to $\alpha(s)\alpha'(s)$ increasing
in $s$, which is true since~$\alpha^2(s)$ is a convex function. This proves (\ref{23dec2735}).

Finally, we prove (\ref{24jan308}). Note that by (\ref{24jan306}) we have 
\be
s_m(I(t,x))\le s_m\Big(\farc{1}{\alpha'(1)}e^{-(x-\eta_\ell(t))}\Big).
\ee
Thus, (\ref{24jan308}) would follow if we show that
\be\label{24jan3020}
s_m\Big(\farc{\xi}{\alpha'(1)}\Big)\le\xi^2,~~\hbox{ for all $0<\xi\le 1$.}
\ee
Applying the decreasing function $\alpha'$ to both sides of (\ref{24jan3020}) and using the definition (\ref{24jan2908}) 
of $s_m(I)$ 
we conclude that (\ref{24jan3020}) is 
equivalent to 
\be\label{24jan3023}
\farc{\alpha'(1)}{\xi}\ge\alpha'(\xi^2),~~\hbox{ for all $0<\xi\le 1$.}
\ee
Recalling that, by assumptions of the present lemma, the function $\beta(s) =\sqrt{s}\alpha'(s)$ is increasing for~$0\le s\le 1$,
we conclude that (\ref{24jan3023}) holds, finishing the proof of (\ref{24jan308}).~$\Box$

\begin{rem}{\rm From now on, we will always assume that  $\alpha(s)$ 
is concave with $\alpha(0)=0$,~$\alpha'(0)=+\infty$,
and, in addition, that the function $\alpha^2(s)$ is convex, and 
the function $\beta(s)=\sqrt{s}\alpha'(s)$ is increasing. 
As we have mentioned, these properties hold for $\alpha(s)=\alpha_1s^k$
with $1/2\le k<1$. 
}
\end{rem}

\subsection{Regularity bounds}\label{sec:regular} 

We collect in this section some a priori regularity bounds on the solutions to the Lucas-Moll system.
\begin{lem}\label{cor-24jan2102}
Suppose that $\alpha(s)$ has the form 
\be\label{24apr2612}
\alpha(s)=\alpha_1s^k,~~\hbox{ with $k\in[1/2,1)$.}
\ee
Then, we have
\be\label{24jan2110}
I_x(t,\eta_\ell(t))\le -\frac{1}{\alpha'(1)}=-\frac{1}{k\alpha_1},
\ee
and  
\be\label{24jan3024}
0\le -I_x(t,x)\le\frac{1}{\alpha'(1)}+\farc{1}{\rho-\kappa}=\farc{1}{k\alpha_1}+\farc{1}{\rho-\kappa},~~\hbox{ for $x>\eta_\ell(t)$}.
\ee
We also have
\be\label{24jan2112}
0\le -s^*_x(t,x) \le K_\alpha:=\farc{1}{1-k}\Big(1+\farc{k\alpha_1}{\rho-\kappa}\Big),~~\hbox{for all $x>\eta_\ell(t)$,}
\ee
and
\be\label{24feb514}
0\le -(\alpha(s^*))_x(t,x) \le k\alpha_1K_\alpha ,~~\hbox{for all $x>\eta_\ell(t)$,}
\ee
\end{lem}
{\bf Proof.}  Note that by assumption (\ref{24apr2612}) about the   form of $\alpha(s)$,
we know that it satisfies all the assumptions of Lemma~\ref{lem-mar502-right}. 
Thus, the bound (\ref{24jan2110}) is an immediate consequence of (\ref{24jan3018bis}) and~(\ref{24jan306}).
Let us also observe that for~$x>\eta_\ell(t)$ we have 
\be
0\le -I_x(t,x)=I(t,x)+\farc{1}{\rho-\kappa}
w(t,x)F(t,x)\le\frac{1}{\alpha'(1)}+\farc{1}{\rho-\kappa}=\farc{1}{k\alpha_1}+\farc{1}{\rho-\kappa},
\ee
which is (\ref{24jan3024}). We used (\ref{24jan2110}) above as well as the monotonicity of $I(t,x)$ and the uniform bounds on $F(t,x)$ and $w(t,x)$. 

Moreover, assumption (\ref{24apr2612}) implies that $s_m(I)$ is given by expression
(\ref{24jan3012}) in Appendix~\ref{sec:alpha-powers} and 
\be
I(t,x)\le 1/(k\alpha_1),~~\hbox{ for $x>\eta_\ell(t)$.}  
\ee
Then, we can estimate 
\be\label{24jan3025}
\bal
s_m'(I)&=\farc{1}{1-k}(k\alpha_1)^{1/(1-k)}I^{k/(1-k)}\le 
\farc{1}{1-k}(k\alpha_1)^{1/(1-k)}\Big(\farc{1}{k\alpha_1}\Big)^{k/(1-k)}\\
&=\farc{k\alpha_1}{1-k},
~~\hbox{for } I\le 1/(k\alpha_1). 
\enbal
\ee
We deduce that
\be\label{24jan3026}
\bal
0&\le -s^*_x(t,x)=-s_m'(I(t,x))I_x(t,x)\le \farc{k\alpha_1}{1-k}
\Big(\frac{1}{\alpha_1k}+\farc{1}{\rho-\kappa}\Big)=\farc{1}{1-k}\Big(1+\farc{k\alpha_1}{\rho-\kappa}\Big), 
\enbal
\ee
for all $x>\eta_\ell(t)$,
which is (\ref{24jan2112}). 

To obtain (\ref{24feb514}), we note that, 
by (\ref{eqn:sm}), \eqref{24jan3025} and expressions \eqref{24apr2606} and \eqref{24feb516} 
in Appendix~\ref{sec:append-deriv}, 
we have 
\be
\bal
0\le-\partial_x\alpha(s^*(t,x))&=-\alpha'(s_m(I(t,x))s_m'(I(t,x))I_x(t,x)=-\farc{s_m'(I(t,x))}{I(t,x)}
I_x(t,x)\\
&=\farc{1}{1-k}(k\alpha_1)^{1/(1-k)}I^{-1+k/(1-k)}(t,x)\Big(I(t,x)+\farc{1}{\rho-\kappa}
w(t,x)F(t,x)\Big)\\
&\le \farc{k\alpha_1}{1-k}+\farc{1}{(1-k)(\rho-\kappa)}(k\alpha_1)^{1/(1-k)+1-k/(1-k)}=
\farc{k\alpha_1}{1-k}\Big(1+\farc{k\alpha_1}{\rho-\kappa}\Big).
\enbal
\ee
We used above the assumption $k\ge 1/2$ so that $k/(1-k)-1\ge 0$.~$\Box$ 
%


\begin{prop}\label{prop-feb102}
Assume that $\alpha(s)$ is given by (\ref{24apr2612})  
and that the initial 
condition~$F_0(x)$ and the terminal condition $w_T(x)$ satisfy (\ref{24jan2322})
and (\ref{24jan3016}), respectively. Suppose also that $0\le\alpha_1\le A$. 
Then, there exists a constant $K_A>0$ that depends on $A$ and $\kappa$ but neither on $T$ nor on $F_0$ or $w_T$, 
such that
\be\label{24feb102}
0\le -F_x(t,x)\le K_A,~~0\le w_x(t,x)\le K_A,~~\hbox{ for all $x\in\Rm$ and $1\le t\le T-1$}.
\ee
\end{prop}
{\bf Proof.}  The evolution equation for the function $F(t,x)$ in (\ref{24sep402}) has the form
\be\label{24feb104} 
\bal
&\pdr{F(t,x)}{t}=\kappa \pdrr{F(t,x)}{x}+c(t,x)F(t,x),
\enbal
\ee
with the function 
\be\label{24feb105}
c(t,x)=\int_{-\infty}^x\alpha(s_m(I(t,y))(-F_y)(t,y)dy,
\ee
that satisfies the a priori bound
\be\label{24feb105b}
0<c(t,x)\le 
\alpha_1 F(t,x)(1-F(t,x)) \le \alpha_1.
\ee
Moreover, we know that 
\be\label{24apr2614}
\hbox{$0< F(t,x)<1$ for all $x\in\Rm$ and $t\in[0,T]$. }
\ee
The Krylov-Safonov bounds~\cite{Krylov-Safonov} 
(see also notes by Sebastien Picard~\cite{Picard}, Theorem 10) 
imply that there exists a universal 
constant~$C_A>0$ and $\delta\in(0,1)$ so that $F(t,x)$ satisfies a H\"older bound
\be\label{24feb106}
\|F\|_{C^{\delta/2,\delta}[1,T]\times\Rm}\le C_A(1+\alpha(1)).
\ee
Next, we integrate by parts to represent the coefficient $c(t,x)$ as 
\be
\bal
&c(t,x)=c_1(t,x)+c_2(t,x),\\
&c_1(t,x)=\alpha(s_m(I(t,x))(1-F(t,x)),\\
&c_2(t,x)=\int_{-\infty}^x[\alpha(s_m(I(t,y))]_y(1-F(t,y))dy.
\enbal
\ee
The function $c_1(t,x)$ satisfies a uniform in time H\"older bound in the $x$-variable
by (\ref{24feb514}) 
and~(\ref{24feb106}).  In addition, the function $c_2(t,x)$ satisfies a uniform in 
time H\"older bound in the $x$-variable by~(\ref{24feb514}). The constants in both of these bounds depend only on $A$. Therefore, the H\"older in space
constant of~$c(t,x)$ is uniformly bounded in $t\in[0,T]$, also uniformly in $0<\alpha_1\le A$. It follows then from~(\ref{24feb106})
that the same can be said of the forcing term $c(t,x)F(t,x)$ in the right side of (\ref{24feb104}).
As its solution~$F(t,x)$ is uniformly bounded by (\ref{24apr2614}), we 
deduce from the Schauder estimates for the forced heat equation that there is a constant $K_A>0$ that depends only on $A$ 
such that
\be\label{24feb108}
0\le -F_x(t,x)\le K_A,~~\hbox{ for all $x\in\Rm$ and $0\le t\le T$}.
\ee
The second bound in (\ref{24feb102}) holds for a similar reason. The function $w(t,x)$
is a solution to the backward in time parabolic equation
\be
\bal
&\pdr{w(t,x)}{t}+\kappa\pdrr{w(t,x)}{x}+2\kappa w_x(t,x)+(\rho-\kappa)\big[1-s_m(I(t,x))-w(t,x)\big]\\
&~~~~~~~~~~~-\alpha(s_m(I(t,x)))w(t,x)F(t,x)=0.
\enbal
\ee
Its solution satisfies $0\le w(t,x)\le 1$ and the coefficients have bounded derivatives
in the $x$-variable by (\ref{24jan2112}), (\ref{24feb514}) and (\ref{24feb108}). This implies the second bound
in (\ref{24feb102}).~$\Box$

\subsection{A bound on the instantaneous learning front speed}\label{sec:bound-speed} 

In this section, we prove an upper bound on the instantaneous speed of the learning front. 
\begin{lem}\label{lem-24jan2106}
Assume that $\alpha(s)$ is given by (\ref{24apr2612})  
and that the initial 
condition~$F_0(x)$ and the terminal condition $w_T(x)$ satisfy (\ref{24jan2322})
and (\ref{24jan3016}), respectively. 
Suppose  that $0\le\alpha_1\le A$. 
There exists~$M_A>0$ that depends on $\rho$, $\kappa$  and $A$ but neither on $T>0$, nor 
on the initial condition  $F_0(x)$ or the terminal condition $w_T(x)$,  
so that the transition point $\eta_\ell(t)$ satisfies 
\be\label{24jan2120}
\dot\eta_\ell(t)\le M_A,~~\hbox{ for all $1\le t\le T-1$.} 
\ee
\end{lem}
{\bf Proof.} 
It will be slightly more convenient for us to 
consider the Lucas-Moll system in terms of the function 
\be\label{24sep421}
v(t,x)=\farc{1}{\rho-\kappa}w(t,x)e^x.
\ee
To the right of the front location~$\eta_\ell(t)$ it takes the form of equation 
(\ref{23dec2727}) in Appendix~\ref{sec:append-deriv} 
\be\label{24jan3110} 
\bal
&\pdr{F(t,x)}{t}=\kappa \pdrr{F(t,x)}{x}+F(t,x)\int_{-\infty}^x\alpha(s_m(I(t,y))(-F_y)(t,y)dy,\\
&\rho v(t,x)= \pdr{v(t,x)}{t}+\kappa\pdrr{v(t,x)}{x}+e^x\big[1-s_m(I(t,x))\big] -\alpha(s_m(I(t,x)))v(t,x)F(t,x).
\enbal
\ee
The location $\eta_\ell(t)$ is determined by relation 
(\ref{24jan3014}) in Appendix~\ref{sec:alpha-powers}:
\be\label{23dec2806bis}
I(t,\eta_\ell(t))=\farc{1}{k\alpha_1},
\ee 
with
\be\label{23dec2808}
I(t,x)= 
e^{-x}\int_x^\infty v(t,y)F(t,y)dy.
\ee
We see from (\ref{23dec2806bis}) that
\be\label{24jan2114}
\dot\eta_\ell(t)=-\farc{I_t(t,\eta_\ell(t))}{I_x(t,\eta_\ell(t)}.
\ee
Using~\eqref{24feb105}  we compute:
\be\label{24feb521}
\bal
&e^{x}\pdr{I(t,x)}{t}=\int_x^\infty [v_tF+vF_t]dy\\
&= 
\int_x^\infty  F[\rho v-e^y(1-s_m(I))-\kappa v_{yy} +\alpha(s_m(I))vF]dy
+\int_x^\infty v\big(\kappa F_{yy}+c F\big)dy\\
 &=\kappa F(t,x)v_x(t,x)-\kappa v(t,x)F_x(t,x)\\
 &+\int_x^\infty F(t,y)\big[\rho v(t,y)-e^y(1-s_m(I(t,y)))+\alpha(s_m(I(t,y)))v(t,y)F(t,y))+c(t,y) v(t,y)].
\enbal
\ee
We may now use the bounds (\ref{24feb102}) in Proposition~\ref{prop-feb102},
together with the uniform bound~(\ref{24feb518}) on~$w(t,x)$ and the relation
(\ref{24sep421}) between $v(t,x)$ and $w(t,x)$,  
to bound the first two terms in the right side: 
\be\label{24jan2116}
0\le \kappa F(t,x)v_x(t,x)-\kappa v(t,x)F_x(t,x)\le K_Ae^x,~~1\le t\le T-1,
\ee
with a constant $K_A$ that depends on $A$. 
Taking also into account the uniform bound~\eqref{24feb105b}, we deduce from 
(\ref{24feb521}) and (\ref{24jan2116}) that $I(t,x)$ satisfies a differential inequality
\be\label{24feb522}
\bal
\pdr{I(t,x)}{t}&\le K_A' 
+K_A'e^{-x}\int_x^\infty F(t,y)v(t,y)dy=K_A'(1+I(t,x)),~~x\ge \eta_\ell(t). 
\enbal
\ee
As $I(t,\eta_\ell(t))$ satisfies (\ref{23dec2806bis}), we obtain from (\ref{24feb522}) that 
\be\label{24feb523}
I_t(t,\eta_\ell(t))\le K_A'\Big(1+\farc{1}{k\alpha_1}\Big).
\ee
Recall also that by Lemma~\ref{cor-24jan2102}
we have
\be\label{24jan2118}
-I_x(t,\eta_\ell(t))\ge \frac{1}{k\alpha_1},
\ee
Taking into account (\ref{24jan2114}), we obtain 
\be
\dot\etal(t)=\farc{I_t(t,\etal(t))}{(-I_x)(t,\etal(t))}\le K_A'(1+k\alpha_1)=M_A,
\ee
and 
(\ref{24jan2120}) follows.~$\Box$

\subsection{Tightness of the learning front}\label{sec:tight-learning} 

Here, we establish a lower bound for $w(t,x)$ to the right of the learning front $\eta_\ell(t)$.  
In other words, the learning front has a width that does not grow in time. 
\begin{prop}\label{lem-jan2202} Under the assumptions of Lemma~\ref{lem-24jan2106}, 
suppose  that $0\le\alpha_1\le A$, with some $A>0$. 
For any $\gamma\in(0,1)$ there exist $T_\gamma>0$, $L_{A,\gamma}>0$ and~$t_\gamma>0$ such that for all $T\ge T_\gamma$ we have 
\be\label{24jan2308}
w(t,x)\ge \gamma,~~\hbox{for all $x>\eta_\ell(t)+L_{A,\gamma}$ and $0\le t\le T-t_\gamma$.}
\ee
These constants do not depend on the terminal condition $w_T(x)$. 
\end{prop}
{\bf Proof.} Let us recall that $w(t,x)$ satisfies  (\ref{24sep402}):
\be\label{24feb526} 
\bal
&w_t+\kappa w_{xx}+2\kappa w_x+(\rho-\kappa)\big[1-s_m(I(t,x))-w\big]
-\alpha(s_m(I(t,x)))wF=0,~~\hbox{ for $x>\eta_\ell(t)$.}
\enbal
\ee
%
Moreover,  
we have, from (\ref{24jan308}) that
\be\label{24jan2204}
s_m(I(t,x))\le e^{-2(x-\eta_\ell(t))},~~
\hbox{ for all $x>\eta_\ell(t)$.}
\ee
Therefore, for any $\delta\in(0,1)$ there exists $L_\delta>0$ that depends solely on $\delta>0$ so that we have
\be\label{24feb527} 
\bal
&w_t+\kappa w_{xx}+2\kappa w_x +
(\rho-\kappa)\big[1- {\delta}-w\big]
-\alpha(\delta)w<0,~\hbox{ $x>\eta_\ell(t)+L_\delta$,}\\
&w(t,\eta_\ell(t)+L_\delta)>0.
\enbal
\ee
Setting 
\be
w(t,x)=\tilde w(t,x-\eta_\ell(t)-L_\delta)
\ee
gives 
\be\label{24feb528} 
\bal
&\tw_t+\kappa \tw_{xx}+2\kappa \tw_x +
(\rho-\kappa)\big[1-\delta-\tw\big]
-\alpha(\delta)\tw<\dot\eta_\ell(t)\tw_x,~\hbox{ $x>0$,}\\
&\tw(t,0)>0.
\enbal
\ee
As $\tw_x>0$ and $\dot\eta_\ell(t)<M_A$ by Lemma~\ref{lem-24jan2106}, 
we have 
\be
\dot\etal(t)\tw_x<M_A\tw_x.
\ee
Thus, the function $\tilde w(t,x)$ satisfies the differential inequality
\be\label{24feb529} 
\bal
&\tw_t+\kappa \tw_{xx}+(2\kappa -M_A)\tw_x +
(\rho-\kappa)\big[1-\delta-\tw\big]
-\alpha(\delta)\tw<0,~\hbox{ $x>0$,}\\
&\tw(t,0)>0.
\enbal
\ee

Let us fix $t_0\in(0,T)$ and $\bar L>0$ and consider a function of the form
\be\label{24jan2206bis}
\bal
z(t,x)&=\psi(t) 
\big(1-e^{-\mu x}\big),~~x>0,~t\le t_0, 
\enbal
\ee
with a decreasing function $\psi(t)$ to be chosen so that $\psi(t_0)=0$ and $\psi(t)>0$ for $t<t_0$. 
Note that~$z(t,x)$ satisfies
\be\label{24jan2208}
\bal
&z(t,0)=0<\tw(t,0),\\
&z(t_0,x)=0< \tw(t_0,x).
\enbal
\ee
We also have
\be\label{24feb530}
\bal
{\cal N}[z]:=z_t+\kappa z_{xx}&+(2\kappa -M_A)z_x +
(\rho-\kappa)\big[1-\delta-z\big]
-\alpha(\delta)z
= \dot\psi(t)(1-e^{-\mu x})-\kappa\mu^2\psi(t)e^{-\mu x}\\
&+(2\kappa-M_A)\psi(t)\mu e^{-\mu x}+
\big(1-{\delta}\big)(\rho-\kappa) -(\rho-\kappa+\alpha(\delta))\psi(t)(1-e^{-\mu x})\\
&=  
\big(1-\delta\big)(\rho-\kappa)-(\rho-\kappa+\alpha(\delta))\psi(t){+ \dot\psi(t)(1-e^{-\mu x})}\\
&+e^{-\mu x}(-\kappa\mu^2+2\kappa\mu-M_A\mu+\rho-\kappa+\alpha(\delta))\psi(t). 
\enbal
\ee
If we take 
\be
\psi(t)=\gamma_1\Big(1-e^{-\gamma_2(t_0-t)}\Big),
\ee
then (\ref{24feb530}) gives 
\be\label{24apr1702}
\bal
{\cal N}[z]&=\big(1-\delta\big)(\rho-\kappa)-(\rho-\kappa+\alpha(\delta))\gamma_1\big(1-e^{-\gamma_2(t_0-t)}\big)
-\gamma_1\gamma_2 e^{-\gamma_2(t_0-t)}(1-e^{-\mu x})\\
&+e^{-\mu x}(\rho-\kappa -\kappa\mu^2+2\kappa\mu-M_A\mu+\alpha(\delta))\gamma_1\big(1-e^{-\gamma_2(t_0-t)}\big)>0,
\enbal
\ee
as long as
\be\label{24jan410}
(1-\delta)(\rho-\kappa)>(\rho-\kappa+\alpha(\delta))\gamma_1+\gamma_1\gamma_2,
\ee
and
\be\label{24apr2616}
\rho-\kappa -\kappa\mu^2+2\kappa\mu-M_A\mu>0.
\ee
Given any $\gamma_1\in(0,1)$, condition  (\ref{24jan410}) holds if we take $\delta>0$ sufficiently small, and then $\gamma_2>0$ also sufficiently small, 
both of them depending on $\gamma_1$. Furthermore, 
(\ref{24apr2616}) holds as long as we take~$\mu>0$ sufficiently small, depending on~$A>0$ but not on $\gamma_1$.

With this choice of $\psi(t)$, both (\ref{24jan2208}) and
(\ref{24apr1702}) hold. In that case, the comparison principle implies 
\be
\tw(t,x)\ge z(t,x),~~\hbox{ for all $0\le t\le t_0$ and $x>0$.}
\ee
That is, we have
\be
\bal
w(t,x)&= \tilde w(t,x-\eta_\ell(t)-L_\delta)\ge z(t,x-\eta_\ell(t)-L_\delta)\\
&=
\gamma_1\Big(1-e^{-\gamma_2(t_0-t)}\Big)\big(1-e^{-\mu(x-\eta_\ell(t)-L_\delta)}\big).
\enbal
\ee
Now,  (\ref{24jan2308}) follows.~$\Box$

\subsection{Tightness of the particles front}\label{sec:tight-particles} 

In addition to the learning front $\eta_\ell(t)$, 
another location of interest is the median front $\eta_m(t)$ of the agents determined by
\be\label{24feb504}
F(t,\eta_m(t))=\farc{1}{2}.
\ee
The distance between the learning front location $\eta_\ell(t)$ and the median front 
$\eta_m(t)$ may grow in time.  
However, as we now show that the level sets of $F(t,x)$ stay together. Let us define the level sets  
\be
F(t,\Gamma_s(t))=s,~~0<s<1.
\ee
We have the following tightness estimate. 
\begin{prop}\label{prop-24apr1802}
Suppose  that $0\le\alpha_1\le A$. 
Given any $0<s_1,s_2<1$ there exists $L_A(s_1,s_2)$ and~$T_A>0$ so that
\be
0<\Gamma_{s_1}(t)-\Gamma_{s_2}(t)\le L_A(s_1,s_2),~~\hbox{ for all $T_A\le t\le T$.}
\ee
\end{prop}
{\bf Proof.} 
It suffices to show that $F(t,x)$ satisfies a differential inequality of the form
\be\label{24apr1808}
F_t\ge\kappa F_{xx}+g_A(F),
\ee
with a nonlinearity $g(u)$ such that 
\be\label{24apr1810}
g_A(0)=g_A(1)=0,~~g_A(u)>0,~~\hbox{ for all $u\in(0,1)$.}
\ee
As the function $F(t,x)$ is a solution to the first equation in (\ref{24jan3110}):
\be\label{24apr1806} 
\bal
&\pdr{F(t,x)}{t}=\kappa \pdrr{F(t,x)}{x}+F(t,x)\int_{-\infty}^x\alpha(s_m(I(t,y))(-F_y)(t,y)dy,
\enbal
\ee
we already know that  for $x<\etal(t)$ the function $F(t,x)$ satisfies the Fisher-KPP equation
\be\label{24apr1814}
F_t=\kappa F_{xx}+\alpha(1)F(1-F).
\ee
In  the region $x>\etal(t)$ the function $F(t,x)$ satisfies
\be\label{24apr1812} 
\bal
&\pdr{F(t,x)}{t}=\kappa \pdrr{F(t,x)}{x}+\alpha(1)F(t,x)(1-F(t,\etal(t)))+F(t,x)\int_{\etal(t)}^x\alpha(s_m(I(t,y))(-F_y)(t,y)dy.
\enbal
\ee
To estimate the integral on the right side of (\ref{24apr1812}),
we
have the following lemma.
\begin{lem}\label{lem-24jan2304} 
For every $p>1$ there exists $k_{p,A}>0$ that depends  on $p$ and $A$, 
and a time~$t_A>0$ that depends also on the other  parameters
of the problem but not on $T$ such that 
\be\label{24jan1910}
\bal
I(t,x)\ge  k_{p,A}F^p(t,x),~~\hbox{ for all $x\ge \etal(t)$ and $t_A\le t\le T$.}
\enbal
\ee
\end{lem}
{\bf Proof.} Equation (\ref{24apr1806}) has the form
\be
F_t=\kappa F_{xx}+c(t,x)F,
\ee
with a function $c(t,x)$ given by~\eqref{24feb105} that satisfies~\eqref{24feb105} and~\eqref{24apr2614}. 
Therefore, according to the same-time Harnack inequality~\cite{BouHenR}, 
for any $1<p<+\infty$, there exist  universal 
constants~$C_{p}>0$ and~$\beta_{p}>0$ that depend only on~$p$, such that for all $t\ge s\ge 1$ and all~$x,y\in\Rm$, we have
\be\label{24jan1912}
F(t,x+y)\ge C_{p}e^{-p\alpha_1(t-s)-p\beta_{p} y^2/[\kappa(t-s)]}F^p(t,x).
\ee
We now optimize over $s\in[1,t]$ for a given $y\in\Rm$ in (\ref{24jan1912}). The function
\be
\phi(s)=\alpha_1s+\farc{\beta_py^2}{\kappa s}
\ee
attains its global minimum at
\be
s(y)=c_0^{-1}|y|,~~c_0^{-1}:=\Big(\farc{\beta_p}{\kappa\alpha_1}\Big)^{1/2}.
\ee
If $|y|<c_0(t-1)$, then we can take $s=t-|y|/c_0$ in (\ref{24jan1912}), 
which gives
\be
F(t,x+y)\ge C_{p} e^{-2pc_0^{-1}\alpha_1|y|}F^p(t,x),~~\hbox{for $|y|\le c_0(t-1)$}.
\ee
Let us now fix some $\gamma\in(0,1)$ and choose $L_{A,\gamma}\ge 0$ as  in
Proposition~\ref{lem-jan2202}, so that (\ref{24jan2308}) holds. 
This gives the following lower bound  for $y>\eta(t)$ and $t>t_{A,\gamma}$:
\be\label{24jan1910bis}
\bal
I(t,y)&=e^{-y}\int_y^\infty e^zw(t,z) F(t,z)dz   \ge 
\gamma\int_{\max(y,\etal(t)+L_{A,\gamma})}^{y+c_0t} e^{z-y}
e^{-2pc_0^{-1}\alpha_1(z-y)}F^p(t,y)dz\\
&\ge k_{p,A}F^p(t,y), 
\enbal
\ee
which is (\ref{24jan1910}).~$\Box$

We go back to the proof of Proposition~\ref{prop-24apr1802}. 
Lemma~\ref{lem-24jan2304} and expression (\ref{24apr2618}) for $\alpha(s_m(I))$ allow us to estimate the integral in the right side of (\ref{24apr1812}) as 
\be
\bal
&\int_{\etal(t)}^x \alpha(s_m(I(t,y))(-F_y)(t,y)dy\ge  \int_{\etal(t)}^x \alpha(s_m(k_{p,A}F^p(t,y)))(-F_y)(t,y)dy\\
&=
c_{p,A}\big({\cG}(F(t,\etal(t))-\cG(F(t,x))\big).
\enbal
\ee
Here, we have set
\be
\cG(u)=\int_0^u\alpha(s_m(z^p))dz.
\ee
We see that for $x>\etal(t)$ the function $F(t,x)$ satisfies a differential inequality
\be\label{23sep2812}
F_t(t,x)\ge \kappa F_{xx}(t,x)+F(t,x)\big[\alpha_1 (1-F(t,\etal(t)))+ c_{p,A} \big({\cG}(F(t,\etal(t))-\cG(F(t,x))\big)\big]. 
\ee
Let us set, for $0<w<v<1$:
\be
q(v)=\alpha_1(1-v)+c_{p,A}\big({\cG}(v)-\cG(w)\big).
\ee
If $c_{p,A}>0$ is sufficiently small, then $q(v)$ is decreasing in $v$ for $v\in(0,1)$. 
Therefore, it attains its minimum over $v\in(w,1)$ at the point $v=1$. It follows that 
\be
\alpha_1 (1-F(t,\etal(t)))+ c_p \big({\cG}(F(t,\etal(t))-\cG(F(t,x))\big)\ge c_p\big({\cG}(1)-\cG(F(t,x))\big),~~\hbox{ for $x>\etal(t)$},
\ee
as long as $c_p>0$ is sufficiently small. 
Using this in (\ref{23sep2812}) gives
\be\label{23sep2814}
\bal
F_t\ge \kappa F_{xx}+c_{p,A}F\big({\cG}(1)-\cG(F(t,x))\big),
~~x>\etal(t). 
\enbal
\ee
In addition, we know that $F(t,x)$ satisfies (\ref{24apr1814}) for $x<\etal(t)$. Together, (\ref{24apr1814}) and (\ref{23sep2814}) give
\be
\bal
&F_t=\kappa F_{xx}+\alpha_1F(1-F),~~x<\etal(t), \\
&F_t\ge\kappa F_{xx}+c_{p,A}F\big({\cG}(1)-\cG(F(t,x))\big) ,
~~x>\etal(t).
\enbal
\ee
Therefore, if we set 
\be\label{24apr1816}
g(s)=s \min\Big(\alpha_1(1-s), c_{p,A}\big({\cG}(1)-\cG(s)\big)\Big),
\ee
then $F(t,x)$ satisfies the differential inequality (\ref{24apr1808}):
\be
F_t\ge\kappa F_{xx}+g(F),~~x\in\Rm. 
\ee
Moreover, the function $g(s)$ defined by (\ref{24apr1816}) satisfies (\ref{24apr1810}). This finishes the proof.~$\Box$
%

\subsection{The intrinsic learning front location}\label{sec:learning-intrinsic}

Let us now describe an approximation to the learning front location in terms of the agent distribution
function $F(t,x)$ alone. Consider the ``intrinsic''  pay-off functional 
\be\label{24feb704}
J(t,x)=J[F](t,x) =\farc{1}{\rho-\kappa}e^{-x}\int_x^\infty e^y F(t,y)dy.
\ee
As $J(t,x)$ is monotonically decreasing, there exists a unique {\it intrinsic learning front location} $e_{\ell}(t)$ such that
\be\label{24feb718}
J(t,e_\ell(t))=\farc{1}{\alpha'(1)}.
\ee
The following proposition shows that $\el(t)$ is a good approximation to the learning front $\etal(t)$.
We will refer to $\el(t)$ as the intrinsic learning front. 
\begin{prop}\label{prop-24feb702}
Suppose  that $0\le\alpha_1\le A$, with some $A>0$. 
There exist $L_{\ell,A}>0$ and $t_A>0$ that do not depend on $T$,so that
\be\label{24feb721}
\el(t)-L_{\ell,A}\le\etal(t)\le \el(t),~~\hbox{ for all $0\le t\le T-t_M$.}
\ee
\end{prop}
{\bf Proof.} 
Comparing to the expression for $I(t,x)$ in (\ref{24sep402})  and recalling that $0\le w(t,x)\le 1$ for all~$x\in\Rm$ and~$0\le t\le T$, we 
immediately see that
\be\label{24feb708}
I(t,x)\le J(t,x),~~\hbox{for all $0\le t\le T$ and $x\in\Rm$.}
\ee
As a consequence, we have the upper bound
\be\label{24feb710}
\eta_\ell(t)\le e_\ell(t),~~\hbox{for all $0\le t\le T$.}
\ee
On the other hand, for any $\gamma\in(0,1)$, we also have, by Proposition~\ref{lem-jan2202} 
\be\label{24feb712}
\bal
I(t,x)&=\farc{1}{\rho-\kappa}e^{-x}\int_x^\infty e^yw(t,y)F(t,y)dy\ge \farc{\gamma}{\rho-\kappa}e^{-x}\int_x^\infty e^yF(t,y)dy\\
&=\gamma J(t,x),~~
\hbox{for all $0\le t\le T-t_\gamma$ and $x\ge\etal(t)+L_{A,\gamma}$.}
\enbal
\ee
Let us fix some $\gamma\in(0,1)$. 
Recalling the upper bound (\ref{24jan306}) on $I(t,x)$ to the right of the learning front, we obtain  from (\ref{24feb712}) that 
\be\label{24feb714}
J(t,x)\le \farc{1}{\gamma\alpha'(1)}e^{-(x-\etal(t))},~~\hbox{for all $0\le t\le T-t_\gamma$ and $x\ge\etal(t)+L_{A,\gamma}$.}
\ee
Thus, if $\el(t)\ge\etal(t)+L_{A,\gamma}$, then 
\be\label{24feb716}
J(t,e_\ell(t))\le\farc{1}{\gamma\alpha'(1)}e^{-(\el(t)-\etal(t))}.
\ee
Comparing to (\ref{24feb718}), we deduce that there exists $L_{\ell,A}$ so that 
\be\label{24feb720}
\el(t)\le \etal(t)+L_{\ell,A},~~ \hbox{for all $0\le t\le T-t_\gamma$,}
\ee
finishing the proof.~$\Box$

One advantage of the intrinsic learning front location $\el(t)$ is that it always moves to the right. 
\begin{prop}\label{lem-24jan1302}
The function $J(t,x)$ defined by (\ref{24feb704}) is monotonically increasing in time and
\be\label{24jan1306}
J(t,x)\ge e^{\kappa(t-s)}J(s,x),~~\hbox{ for all $0<s<t$.}
\ee
Moreover, the point $\el(t)$ always moves to the right: 
\be\label{24jan1314}
\dot\el(t)\ge\kappa.
\ee
\end{prop}
{\bf Proof.} The function $F(t,x)$ is a super-solution to the heat equation
\be\label{24jan1302}
F_t\ge\kappa F_{xx},
\ee
and is non-increasing.
 Let us multiply (\ref{24jan1302}) by $e^x$ and integrate from $x$ to $+\infty$. This gives
\be
\bal
(\rho-\kappa)e^x\pdr{J(t,x)}{t}&\ge\kappa\int_x^\infty e^yF_{yy}(t,y)dy=-\kappa\int_x^\infty e^yF_y(t,y)dy-\kappa e^xF_x(t,x)\\
&=\kappa\int_x^\infty e^yF(t,y)dy+\kappa e^x  F(t,x)-\kappa e^xF_x(t,x)\ge (\rho-\kappa)\kappa e^xJ(t,x).
\enbal
\ee
This implies  (\ref{24jan1306}). It also follows that 
\be\label{24apr1908}
\pdr{J(t,\el(t))}{t}\ge\kappa J(t,\el(t))+\farc{\kappa}{\rho-\kappa}F(t,\el(t))=\farc{\kappa}{\alpha'(1)}+\farc{\kappa}{\rho-\kappa}F(t,\el(t)).
\ee
In addition, as 
\be
-J_x(t,x)=J(t,x)+\farc{1}{\rho-\kappa}F(t,x),
\ee
we know that
\be\label{24apr1910}
-J_x(t,\el(t))= J(t,\el(t))+\farc{\kappa}{\rho-\kappa}F(t,\el(t))=\farc{1}{\alpha'(1)}+\farc{\kappa}{\rho-\kappa}F(t,\el(t)).
\ee
It follows from (\ref{24feb718}), (\ref{24apr1908}) and (\ref{24apr1910})  that 
\be
\dot \el(t)=-\farc{J_t(t,\el(t))}{J_x(t,\el(t)}\ge\kappa,
\ee
finishing the proof ~$\Box$

A corollary of Propositions~\ref{prop-24feb702} and~\ref{lem-24jan1302} is the following lower bound on the learning front location. 
\begin{cor}\label{cor-24jan1302}
There exists a time $T_0$ 
and two constants $M_A>0$ and $t_A>0$ so that for all~$T\ge T_0$  we have
\be\label{24jan2304}
\eta_\ell(t)\ge\kappa t-M_A,~~\hbox{ for all $0\le t\le T-t_A$.}
\ee
\end{cor}

One may ask how the median front location $\eta_m(t)$ defined by (\ref{24feb504}):
\be\label{24feb722}
F(t,\eta_m(t))=\farc{1}{2}.
\ee
compares to the intrinsic  learning front location $\el(t)$ and the learning front location $\etal(t)$.  
 We have, for any $x>\el(t)$:
\be\label{24mar325}
\bal
\farc{1}{\alpha'(1)}&=J(t,\el(t))=\farc{e^{-\el(t)}}{\rho-\kappa}\int_{\el(t)}^\infty e^yF(t,y)dy\ge \farc{e^{-\el(t)}}{\rho-\kappa}F(t,x)
\int_{\el(t)}^x e^ydy\\
&=F(t,x)\frac{1}{\rho-\kappa}\Big(e^{x-\el(t)}-1\Big).
\enbal
\ee
As $F(t,\eta_m(t))=1/2$, we deduce from (\ref{24mar325}) and Proposition~\ref{prop-24feb702} the following,
\begin{prop}\label{prop-24apr1904}
There exist $T_0$, $L_{A,\ell}>0$, and $t_A>0$ so that for all~$T\ge T_0$  we have
\be\label{24mar324}
\eta_m(t)\le\etal(t)+L_{A,\ell},~~~~\hbox{ for all $0\le t\le T-t_A$.}
\ee
\end{prop}

\section{Spreading in the lottery society}\label{sec:lottery}

In this section we prove Theorem~\ref{conj-feb2202_intro}. 

\subsection{The lottery society}\label{sec:lottery-intro}

The lottery society is the regime when the median front is located far behind the learning front:
\be\label{24apr1918}
\eta_m(t)\ll \eta_\ell(t),
\ee  
and
\be\label{24jun1102}
\eta_l(t)-\eta_m(t)\to+\infty,~\hbox{as $t\to+\infty$.}
\ee
In particular, it follows that 
\be\label{24apr1920}
s_m(t,x)= 1 \hbox{ unless $F(t,x)\ll 1$.}
\ee
In other words, in the lottery regime the overwhelming majority of the agents are only searching
and are not producing anything at all. Moreover, the small fraction of those who do produce is diminishing in time
because the learning and median fronts diverge, as seen in (\ref{24jun1102}).

We will be more precise below but for now let us informally analyze when one can expect the lottery society to form. 
If (\ref{24apr1920}) holds, then $F(t,x)$ solves, approximately, the Fisher-KPP equation 
\be\label{24apr1921}
F_t\approx \kappa   F_{xx}+\alpha_1F(1-F).
\ee
Therefore, the median front has the asymptotics
\be\label{24feb2925}
\eta_m(t)=c_*t+o(t),~~c_*=2\sqrt{\kappa\alpha_1}.
\ee
Neglecting the algebraic and logarithmic terms in $t$ and $x$, 
the solution to (\ref{24apr1921}) has the exponential asymptotics 
\be
F(t,x+c_*t)\sim e^{-\lambda _*x}e^{-|x|^2/(4\kappa t)},~x\gg c_*t,
\ee
with
\be\label{24sep902}
\lambda_*=\sqrt{\farc{\alpha_1}{\kappa}}=\farc{c_*}{2\kappa}.
\ee
For reasons that will become clear very soon, we
expect the lottery society to exist as long as the exponential decay rate satisfies
\be\label{24feb1902}
\lambda_*<1,
\ee
that is, as long as 
\be
c_*<2\kappa,
\ee
or, equivalently,
\be\label{24feb1908}
\alpha_1<\kappa.
\ee
Let us now informally compute the asymptotics of the learning front location when
the exponential decay rate does satisfy (\ref{24feb1902}). In that regime, the 
intrinsic pay-off functional $J(t,x)$ from~\eqref{24feb704} has 
the asymptotics, for $x\gg \eta_m(t)\approx c_*t$:
\be
\bal
J(t,x)&= \farc{1}{\rho-\kappa}e^{-x}\int_x^\infty e^yF(t,y)dy \approx \farc{1}{\rho-\kappa}e^{-x}\int_x^\infty e^y
e^{-\lambda _*(y-c_*t)}e^{-|y-c_*t|^2/(4\kappa t)}dy\\
&= \farc{1}{\rho-\kappa}e^{-x}\int_x^\infty e^y
e^{-\lambda _*(y-c_*t)}e^{-|y-c_*t|^2/(4\kappa t)}dy\\
&=\farc{1}{\rho-\kappa}e^{-x+c_*t}\int_{x-c_*t}^\infty e^{(1-\lambda_*)y}e^{-|y|^2/(4\kappa t)}dy.
\enbal
\ee
Assuming that (\ref{24feb1902}) holds, once again,  dropping algebraic factors in $t$ and constants, gives 
\be
\bal
J(t,x)&\sim e^{-x+c_*t}\int_{(x-c_*t)/\sqrt{4\kappa t}}^\infty e^{(1-\lambda_*)\sqrt{4\kappa t}y}e^{-|y|^2 }dy
\\
&= e^{-x+c_*t+(1-\lambda_*)^2\kappa t}\int_{(x-c_*t)/\sqrt{4\kappa t}}^\infty e^{-|y-(1-\lambda_*)\sqrt{\kappa t}|^2 }dy\\
&= e^{-x+c_*t+(1-\lambda_*)^2\kappa t}\int_{\Gamma(t)}^\infty e^{-|y|^2 }dy,
\enbal
\ee
with
\be
\Gamma(t)=\farc{(x-c_*t)}{\sqrt{4\kappa t}}-(1-\lambda_*)\sqrt{\kappa t}.
\ee
Therefore, we have $J(t,x)\gg 1$ if we have both 
\be\label{24feb1904}
-x+c_*t+(1-\lambda_*)^2\kappa t\gg 1,
\ee
and $\Gamma(t)\le 0$, which is
\be\label{24feb1906}
x-c_*t<2(1-\lambda_*)\kappa t
\ee
Note that (\ref{24feb1906}) is a consequence of (\ref{24feb1904}) when $\lambda_*<1$. Thus, one expects 
that if (\ref{24feb1902}) holds then we have
\be
\eta_\ell(t)\sim v_*t,
\ee
with
\be
v_*=c_*+(1-\lambda_*)^2\kappa =2\sqrt{\kappa\alpha_1}+(\sqrt{\kappa}-\sqrt{\alpha_1})^2
=\kappa+\alpha_1.
\ee
According to (\ref{24feb1908}), we should expect that
the transition out of the lottery society happens exactly when~$\alpha_1=\kappa$
and~$v_*=2\kappa$.  This exactly 
matches the lower bound for the existence of balanced growth paths in~\cite{Porretta-Rossi}. 

We now restate  our main result on existence of the lottery society.  
\begin{thm}\label{conj-feb2202}
Consider solutions of~\eqref{eqn:intro1} with  $\alpha_1< \kappa$. 
Suppose the initial condition $F_0(x)$ is monotonically decreasing, 
satisfies $F_0(x)=1$ for all $x\le -L_0$, and $F_0(x)=0$ for all $x\ge L_0$. Suppose 
 the terminal condition $w_T(x)$ is increasing and satisfies $w_T(-\infty)=0$, $w_T(+\infty)=1$.
Then, we have the following asymptotics:
\be\label{24feb2618}
\eta_m(t)=2\sqrt{\kappa\alpha_1}t+o(t),
\ee
and
\be\label{24feb2616}
\etal(t)=(\kappa+\alpha_1)t+o(t).
\ee
\end{thm}
As the Lucas-Moll system is only formulated on a finite time interval $0\le t\le T$, the $o(t)$ terms in Theorem~\ref{conj-feb2202} should be understood as follows:
we say that $g(t)=o(t)$ if  for any $\eps>0$ there exist~$t_0^\eps>0$ and $t_1^\eps>0$ so that for all $T>t_1^\eps+t_0^\eps$ we have 
\be
|g(t)|<\eps,~~\hbox{ for all $0<t_0^\eps<t<T-t_1^\eps$.}
\ee
The rest of this section contains the proof of Theorem~\ref{conj-feb2202}.  
We first obtain the learning front location estimate~\eqref{24feb2616} in Section~\ref{sec:learn}. 
The agents front location estimate~\eqref{24feb2618} is deduced 
from these estimates in Section~\ref{sec:agent-loc}. Section~\ref{sec:lem-aux} contains the proof
of an auxiliary technical result, namely Lemma~\ref{lem-mar304},  used in the estimates on the learning front location. 
 
 \subsection{Proof of the learning front location estimate~\eqref{24feb2616}}\label{sec:learn}
 
 By Proposition~\ref{prop-24feb702} the estimate \eqref{24feb2616} will follow
 from the similar estimate for $\el(t)$:
 \be\label{24feb2616_bis}
\el(t)=(\kappa+\alpha_1)t+o(t).
\ee
In order to obtain the upper and the lower bounds in~\eqref{24feb2616_bis}   we construct super- and sub-solutions of
$F(t,x)$, respectively. 

{\bf Upper bound in~\eqref{24feb2616_bis}.} The upper bound comes 
from the function $\bG(t,x)$,  the solution to the linearized Fisher-KPP equation 
\be\label{24jun1202}
\bal
&\bG_t= \kappa \bG_{xx}+\alpha_1\bG,\\
&\bG(0,x)=F_0(x).
\enbal
\ee
Positivity of the functional $J: F \to \mathbb{R}$ defined in~\eqref{24feb704} 
and the inequality $F(t,x)\le \bG(t,x)$ imply that  $J[F](t,x)\le J[\bG](t,x)$.
Therefore
 \be\label{24feb2621}
\el(t)\le \el^+(t),~~\hbox{ for all $t>0$,}
\ee
where we set $\el^+(t)$ via
\be
J[\bG](t,\el^+(t))=\farc{1}{\alpha'(1)}=\farc{2}{\alpha_1}.
\ee
It remains to show
\be\label{24feb2924}
\el^+(t)\le (\kappa+\alpha_1)t+o(t).
\ee
 Observe that the function 
\be
\bar J(t,x)=J[\bar G](t,x)=\farc{e^{-x}}{\rho-\kappa}\int_x^\infty e^y\bar G(t,y)dy
\ee
satisfies 
\be\label{24jun1110}
\bal
\bar J_t(t,x)-\bar J_{xx}(t,x)&=\farc{e^{-x}}{\rho-\kappa}\int_x^\infty e^y\big[\bar G_{yy}(t,y)+\alpha_1\bar G(t,y)]dy-\bar J(t,x)
-\farc{2e^{-x}}{\rho-\kappa}e^x \bar G(t,x)\\
&+\farc{e^{-x}}{\rho-\kappa}\partial_x(e^x\bar G(t,x))
=\farc{e^{-x}}{\rho-\kappa}\Big[-e^x\bar G_x-\int_x^\infty e^y\bar G_y(t,y)dy\Big]+\alpha_1\bar J(t,x)\\
&-\bar J(t,x)-\farc{2}{\rho-\kappa}\bar G(t,x)+\farc{1}{\rho-\kappa}\bar G(t,x)+\farc{1}{\rho-\kappa}\bar G_x(t,x)= \alpha_1\bar J(t,x).
\enbal
\ee
Moreover, as the function $F_0(x)$ is compactly supported on the right and satisfies $0\le  F_0(x)\le 1$ for all $x\in\Rm$, at the time $t=0$ we have
\be\label{24jun1112}
\bar J(0,x)=\farc{e^{-x}}{\rho-\kappa}\int_x^\infty e^y F_0(y)dy\le \farc{C_0}{\rho-\kappa}e^{-x}.
\ee
It follows from (\ref{24jun1110}) and (\ref{24jun1112}) that
\be
\bar J(t,x)\le \farc{C_0}{\rho-\kappa}e^{-x+(\kappa+\alpha_1)t},
\ee
and (\ref{24feb2924}) follows. This completes the proof of the upper bound in~\eqref{24feb2616_bis}.

{\bf Lower bound in~\eqref{24feb2616_bis}.}  The proof of the lower bound
in~\eqref{24feb2616_bis} is a bit longer.  
We will, in particular, use the upper bound 
\be\label{24feb2908}
F(t,x)\le \bF(t,x).
\ee
Here, $\bF(t,x)$ is the solution to the Fisher-KPP equation 
\be\label{24feb2620}
\bal
&\bF_t= \kappa \bF_{xx}+\alpha_1\bF(1-\bF),\\
&\bF(0,x)=F_0(x).
\enbal
\ee
As the function~$\alpha(s_m(I(t,x)))$
is decreasing in $x$,   we may use  (\ref{24feb2908}) to estimate  
the integral in the equation for $F(t,x)$ from below, as
\be\label{24feb2906}
\bal
\int_{-\infty}^x\alpha(s_m(I(t,y)))&(-F_y)(t,y)dy =
\int_{-\infty}^x\alpha(s_m(I(t,y)))(1-F)_y(t,y)dy\\
&=\alpha(s_m(I(t,x)))(1-F)(t,x)+\int_{-\infty}^x(1-F(t,y))[-\alpha(s_m(I(t,y))]_ydy
\\
&\ge  \alpha(s_m(I(t,x)))(1-F)(t,x)+\int_{-\infty}^x(1-\bar F(t,y))[-\alpha(s_m(I(t,y))]_ydy.
\enbal
\ee
We integrate by parts once again to get
\be\label{24feb2910}
\bal
\int_{-\infty}^x\alpha(s_m(I(t,y)))&(-F_y)(t,y)dy 
\ge  \alpha(s_m(I(t,x)))(1-F)(t,x)-\alpha(s_m(I(t,x)))(1-\bar F)(t,x)\\
&+\int_{-\infty}^x\alpha(s_m(I(t,y))(-\bar F(t,y))_ydy=
\alpha(s_m(I(t,x)))(\bar F-F)(t,x)\\
&+\int_{-\infty}^x\alpha(s_m(I(t,y))(-\bar F(t,y))_ydy\ge 
\int_{-\infty}^x\alpha(s_m(I(t,y))(-\bar F(t,y))_ydy.
\enbal
\ee
Therefore, the function $F(t,x)$ satisfies a differential inequality
\be\label{24feb2912}
F_t\ge \kappa F_{xx}+F(t,x)\int_{-\infty}^x\alpha(s_m(I(t,y))(-\bar F(t,y))_ydy.
\ee
This will be our starting point in getting a lower bound for $\el(t;F)$ in~\eqref{24feb2616_bis}. 
The next step is to replace $I(t,x)$ in the right side of (\ref{24feb2912})
that depends on the   function 
$w(t,x)$ by a functional of $J(t,x)$ that depends solely on $F(t,x)$ and not  on $w(t,x)$. 
To do this, we
first note that 
\be\label{24feb2914}
\int_{-\infty}^x\alpha(s_m(I(t,y))(-\bar F(t,y))_ydy=\alpha_1(1-\bar F(t,x)),~~\hbox{ for $x\le\etal(t)$},
\ee
and
\be\label{24feb2916}
\int_{-\infty}^x\alpha(s_m(I(t,y))(-\bar F(t,y))_ydy\ge
\alpha_1(1-\bar F(t,\etal(t))),
~~\hbox{ for $x\ge\etal(t)$}.
\ee
We also know from Proposition~\ref{prop-24feb702} that
\be\label{24feb2923}
\tel(t;F):=\el(t;F)-L_{\ell,A} \le \etal(t)\le \el(t;F).
\ee
Therefore, we have
\be\label{24feb2918}
\int_{-\infty}^x\alpha(s_m(I(t,y))(-\bar F(t,y))_ydy=\alpha_1(1-\bar F(t,x)),
~~\hbox{ for $x\le\tel(t;F)$}.
\ee
As the function $\bar F(t,x)$ is decreasing, we also know that 
\be
\label{24feb2918bis}
\bal
\int_{-\infty}^x\alpha(s_m(I(t,y))(-\bar F(t,y))_ydy&=\alpha_1(1-\bar F(t,x))\\
&\ge 
\alpha_1(1-\bar F(t,\tel(t;F))),
~~\hbox{ for $\tel(t;F)\le x\le\eta_\ell(t)$},
\enbal
\ee
and
\be\label{24feb2920}
\bal
\int_{-\infty}^x\alpha(s_m(I(t,y))(-\bar F(t,y))_ydy\ge \alpha_1(1-\bar F(t,\eta_\ell(t)))\ge 
\alpha_1(1-\bar F(t,\tel(t;F))),
\hbox{ for $x\ge\etal(t)$}. 
\enbal
\ee
We can summarize the above bounds as
\be\label{24feb2922}
\bal
&\int_{-\infty}^x\alpha(s_m(I(t,y))(-\bar F(t,y))_ydy\ge \alpha_1\bcR[\tel(t;F)](t,x).
\enbal
\ee
Here, for $y\in\Rm$, we have set
\be
\bal
\bcR[y](t,x):=(1-\bar F(t,x))\one(x\le y)+ (1-\bar F(t,y))\one(x\ge y).
\enbal
\ee  
Note that the function $\bcR[y](t,x)$ depends on  $\bar F(t,x)$.
Since the function $\bar F(t,x)$ is decreasing in $x$, we have the monotonicity property
\be\label{24mar102}
\bcR[y_1](t,x)\ge \bcR[y_2](t,x),~~\hbox{ for all $x\in\Rm$ if $y_1>y_2$.}
\ee
In addition, if 
\be
0\le F_1(t,x)\le F_2(t,x)\le 1,~~\hbox{ for all $t>0$ and $x\in\Rm$,}
\ee
then we have 
\be
\tel[t,F_1]\le\tel[t,F_2].
\ee 
Thus,~the function $\bcR[\tel(t;F)](t,x)$ is monotonic in $F$:
\be
\bcR[\tel(t;F_1)](t,x)\le \bcR[\tel(t;F_2)](t,x),\hbox{ for all  $x\in\Rm$ if 
$0\le F_1(x)\le F_2(x)\le 1$ for all  $x\in\Rm$.}
\ee
This gives us a comparison principle.
\begin{lem}\label{lem-24feb2902}
Suppose that $P(0,x)\ge Q(0,x)$ for all $x\in\Rm$, and that $P(t,x)$ and $Q(t, x)$ satisfy 
the differential inequalities 
\be
P_t\ge \kappa P_{xx}+\alpha_1\bcR[\tel(t;P)] P,
\ee 
and
\be
Q_t\le \kappa Q_{xx}+\alpha_1\bcR[\tel(t;Q)]Q.
\ee 
Then, we have
\be
P(t,x)\ge Q(t,x),~~\hbox{ for all $t>0$ and $x\in\Rm$.}
\ee
\end{lem}
On the other hand, we see from (\ref{24feb2912}) and (\ref{24feb2922}) 
that $F(t,x)$ satisfies the differential inequality 
\be\label{24mar106}
F_t\ge \kappa F_{xx}+\alpha_1\bcR[\tel(t;F)]F. 
\ee  
Thus, we have the following.
\begin{cor} \label{cor-24jun1302}
Suppose that $Q(t,x)$ satisfies 
\be\label{24mar104}
\bal
&Q_t\le \kappa Q_{xx}+\alpha_{1}\bcR[\tel(t;Q)]Q,~~t>t_0,\\
&Q(t_0,x)\le F(t_0,x),
\enbal 
\ee  
then $Q(t,x)\le F(t,x)$ for all $t>t_0$.
\end{cor}
Our goal is to construct a solution to the differential inequality (\ref{24mar104}) and estimate the location~$\tel(t;Q)$ from below.
This will give a lower bound on $\tel(t;F)$, which, in turn, will provide a lower bound on $\etal(t)$, via (\ref{24feb2923}).  

To this end, let us recall that $c_*$ is given by (\ref{24feb2925}), set
\be\label{24feb2930}
\bar z(t,x)=\bcR[(c_*+\delta)t-L_1](t,x),
\ee
and let $Q(t,x)$ be a solution to
\be\label{24feb2927}
\bal
&Q_t=\kappa Q_{xx}+\alpha_1\bar z(t,x)Q,~~t>t_0.
\enbal
\ee
Here, $\delta>0$ is chosen so that
\be\label{24feb2926}
c_*+2\delta< \alpha_1+\kappa.
\ee
This is possible because $\alpha_1<\kappa$ and $c_*=2\sqrt{\alpha_1\kappa}$. The right side in (\ref{24feb2926}) is the approximate location  of
the learning front of the (linearized) 
Fisher-KPP solution $\el^+(t)$, as seen from~(\ref{24feb2924}). The shift~$L_1$ in~(\ref{24feb2930}) will be chosen 
later, to accommodate the initial time interval $0\le t\le t_0$.  

We first claim the following. 
\begin{lem}\label{lem-mar304}
There is a solution $Q(t,x)$ to (\ref{24feb2930})-(\ref{24feb2927}) such that $0\le Q(t_0,x)\le F(t_0,x)$ and
\be\label{24mar108}
\tel(t;Q)\ge (c_*+\delta)t-L_1,~~\hbox{ for all $t>t_0$,}
\ee
and
\be\label{24feb2929}
\tel(t;Q)\ge (\alpha_1+\kappa-\delta)t+o(t).
\ee
\end{lem}
Note that (\ref{24mar108}), in turn, implies that 
\be\label{24feb2928}
\bar z(t,x)\le \bcR[\tel(t;Q)](t,x),
\ee
and gives the next corollary.
\begin{cor}\label{cor-24mar306}
The solution  $Q(t,x)$ to (\ref{24feb2930})-(\ref{24feb2927}) with the initial condition
$Q(t_0,x)=F(t_0,x)$ satisfies (\ref{24mar104}). 
\end{cor}
Corollary~\ref{cor-24jun1302} and (\ref{24feb2929}) show then that 
\be\label{24mar316}
\tel(t;F)\ge\tel(t;Q)\ge (\alpha_1+\kappa-\delta)t+o(t),
\ee
proving the lower bound in (\ref{24feb2616_bis}). The proof of Lemma~\ref{lem-mar304} is presented
in Section~\ref{sec:lem-aux}.

 \subsection{Proof of the  agents front location estimate~\eqref{24feb2618}}\label{sec:agent-loc}
  
We now look at the long time asymptotics of the agents front location in Theorem~\ref{conj-feb2202}. 
This, however, is a simple consequence of the already proved lower bound on the learning front location in 
(\ref{24feb2616}). Indeed, it follows from (\ref{24feb2616}) that for any $\delta>0$ there exist $t_0>0$ and $L_0>0$ 
so that 
\be\label{24jun1314}
\etal(t)\ge (\kappa+\alpha_1-\delta)t-L_0,~~\hbox{for all $t>t_0$.}
\ee
Note that none of these parameters depend on the terminal time $T$. Therefore, for all $t>t_0$ the function
$F(t,x)$ satisfies the Fisher-KPP equation
\be
F_t=\kappa F_{xx}+\alpha_1F(1-F),~~x<(\kappa+\alpha_1-\delta)t-L_0,~~\hbox{for all $t>t_0$.}
\ee
In addition, it is well-known  from the standard Fisher-KPP theory that for any $c<c_*=2\sqrt{\kappa\alpha_1}$ there exists a decreasing 
function $\un U(x)$, defined for~$x<0$, that satisfies
\be
0<U(x)<1,~~\hbox{ for all $x<0$}, \lim_{x\to-\ifnty}U(x)=1,
\ee
and
\be
\bal
&-c\un U_x=\kappa \un U_{xx}+\alpha_1\un U(1-\un U),~~x<0,
&\un U(0)=0,
\enbal
\ee
The comparison principle implies that if we choose  $\tilde L_0$ so that 
\be
F(t_0,x)\ge U(x+\tilde L_0),~~\hbox{ for all $x<-\tilde L_0$},
\ee
then 
\be\label{24jun1318}
F(t,x)\ge U(x-ct+\tilde L_0),~~\hbox{ for all $t\ge t_0$ and $x<-\tilde L_0+ct$.}
\ee
As $U(x)$ approaches its limit as $x\to-\infty$ no faster than exponentially,
existence of such $L_0$ follows from the assumption that there exists $\ell_0$ so that $F_0(x)=1$ for all $x\le\ell_0$
and the differential inequality
\be
F_t\ge\kappa F_{xx}.
\ee
We immediately deduce from (\ref{24jun1318}) that 
\be
\eta_m(t;F)\ge ct-\tilde L_0,~~\hbox{ for all $t>t_0$.}
\ee
As  $c<c_*=2\sqrt{\kappa\alpha_1}$ is arbitrary, the lower bound in (\ref{24feb2618}) follows immediately:
\be
\eta_m(t)\ge 2\sqrt{\kappa\alpha_1}t+o(t).
\ee
On the other hand, the upper bound 
\be
\eta_m(t)\le 2\sqrt{\kappa\alpha_1}t+o(t),
\ee
follows immediately from the comparison (\ref{24feb2908}) of $F(t,x)$ to the solution to the Fisher-KPP equation~(\ref{24feb2620}).~$\Box$

\subsection{The proof of Lemma~\ref{lem-mar304}} \label{sec:lem-aux}

\subsubsection{Construction of a sub-solution}

In order to prove Lemma~\ref{lem-mar304} we will construct a sub-solution to to (\ref{24feb2930})-(\ref{24feb2927}).  
Let $Q(t,x)$ be a solution to (\ref{24feb2930})-(\ref{24feb2927}). We first do some standard preliminary transformations: first, set 
\be
Q(t,x)=P(t,x-c_*t)e^{-\lambda_*(x-c_*t)},
\ee
with  $\lambda_*=c_*/(2\kappa)=\sqrt{\alpha_1/\kappa}$
which gives
\be\label{24feb2931}
P_t-c_*P_x+c_*\lambda_* P=\kappa P_{xx}-2\kappa\lambda_*P_x+\kappa\lambda_*^2P+\alpha_1\bar z(t,x+c_*t)P.
\ee
As $c_*=2\kappa\lambda_*$ and 
\be
c_*\lambda_*=\kappa\lambda_*^2+\alpha_1,
\ee
we get from (\ref{24feb2931}) that 
\be\label{24feb2932}
P_t=\kappa P_{xx}-\alpha_1[1-\bar z(t,x+c_*t)]P.
\ee
Observe that the zero-order term above has the form
\be\label{24jun1306}
\bal
1-\bar z(t,x+c_*t)&=
1-(1-\bar F(t,x+c_*t))\one(x+c_*t\le (c_*+\delta)t-L_1)\\
&- (1-\bar F(t,(c_*+\delta)t-L_1))\one(x+c_*t\ge (c_*+\delta)t-L_1)\\
&=\bar F(t,x+c_*t)\one(x\le \delta t-L_1)+
\bar F(t,(c_*+\delta)t-L_1)\one(x\ge \delta t-L_1).
\enbal
\ee
In particular, we have 
\be
0\le 1-\bar z(t,x+c_*t)\le 1,~~\hbox{for all $t>0$ and $x\in\Rm$}.
\ee
Moreover, it follows from the standard estimates for the solution to the Fisher-KPP equation and~(\ref{24jun1306}) that 
there is $\omega>0$ so that for any $\gamma>0$ we have
\be\label{24mar112}
0\le 1-\bar z(t,x+c_*t)\le C_\gamma\exp(-\omega t^\gamma),~~\hbox{for all $x>t^\gamma$.}
\ee
Thus, for $x\gg 1$, a good approximation to the solution of (\ref{24feb2932})
is by a solution to the standard heat equation.

We are going to construct a sub-solution to \eqref{24feb2932} in the form
\be
v(t,x)=\max[0,\eps\underline v(t,x)],~~t>t_0,~~x>0,
\ee
with
\be\label{24mar122}
\un v(t,x)=v_1(t,x)-v_2(t,x),
\ee
and
\be
\bal
&v_1(t,x)=\xi(t)\un u(t,x),\\
&v_2(t,x)=\farc{1}{t^{3/2-\beta}}
\cos\Big(\farc{x}{t^\mu}\Big)\one\big(x\le \farc{3\pi}{2}t^\mu\big).
\enbal 
\ee
Here, $\un u(t,x)$ is a solution to the heat equation on half-line:
\be\label{24mar302}
\bal
&\un u_t=\kappa\un u_{xx},~~x>0, ~t>0,\\
&\un u(t,0)=0, \\
&\un u(0,x)=\phi(x),~~x>0.
\enbal
\ee
The function $\xi(t)$, the parameters $t_0>0$, $\mu\in(1/2,1/3)$, $\eps>0$ and $\beta>0$ will be chosen below,
as will be the initial condition $\phi(x)$ for (\ref{24mar302}).  In particular, we will choose $\xi(t)$ so that 
\be\label{24mar123}
\frac{\xi(t_0)}{2}\le\xi(t)\le\xi(t_0),~~\hbox{for all $t>t_0$},
\ee

We will need to verify two properties: first, that the initial comparison holds at the time $t=t_0$:
\be\label{24mar318}
P(t_0,x)=F(t_0,x)e^{\lambda_*(x-c_*t_0)}\ge  \max[0,\eps\un v(t_0,x)],~~\hbox{for all $x>0$.}
\ee
Second, we will need to check that  
\be\label{24mar320}
v_t\le \kappa v_{xx}-\alpha_1[1-\bar z(t,x+c_*t)]v,~~\hbox{ wherever $v(t,x)>0$ and $x>0$.}
\ee
We will not need to check the boundary condition at $x=0$ because we will see that by construction we will have 
\be\label{24mar321}
v(t,x)=0,~~\hbox{for $x>0$ sufficiently small.}
\ee

We first explain how (\ref{24mar321}) comes about. 
The function
$\un u(t,x)$  behaves as 
\be\label{24mar121}
\un u(t,x)\sim \farc{x}{t^{3/2}},~~\hbox{ for $t>t_0$, $x\ll t^{1/2}$, }
\ee
and is increasing in that region, provided that $t_0$ is sufficiently large and the initial condition
$\phi(x)$ in~(\ref{24mar302}) is chosen appropriately.
In particular, this monotonicity property and asymptotics hold in the region~$x\le (3/2)\pi t^\mu$, where both terms in the right side of~(\ref{24mar122}) are non-zero,  
as long as we take $\mu<1/2$.
On the other hand, the function $v_2(t,x)$ is decreasing in $x$ where it is positive and 
\be
v_2(t,x)\sim \farc{1}{t^{3/2-\beta}}\gg \farc{x}{t^{3/2}},~~\hbox{ for $x\ll t^\beta$.}
\ee
Thus, assuming that (\ref{24mar123}) holds, if we choose~$\xi(t_0)$ sufficiently small (but independently of $t_0$), 
then there is a point~$L_1(t)$ such that
\be\label{24mar116}
t^{\beta/2}<L_1(t)\le \farc{\pi}{2}t^\mu,
\ee
and 
\be
\bal
&\un v(t,L_1(t))=0,\\
&\un v(t,x)\le 0,~~\hbox{ for all~$0<x<L_1(t)$},\\
&\un v(t,x)>0,~~\hbox{for $x>L_1(t)$.}
\enbal
\ee
Therefore, not only (\ref{24mar321}) holds but also
$v(t,x)=0$ in the region $I_L=\{0\le x<L_1(t)\}$, so that, in particular, $v(t,x)$ is a sub-solution for the equation (\ref{24feb2932}) for $P(t,x)$ 
in that region.

 We will consider now separately the remaining the right region 
 \[
 I_R=\{x>(3/2)\pi t^\mu\},
 \]
 and the middle region
\[
I_M=\{L_1(t)<x<(3/2)\pi t^\mu\}.
\] 
 
{\bf Step 1. The  right region.} 
First, in the region $I_R$, we have $v(t,x)=v_1(t,x)=\xi(t)\un u(t,x)$.
Using (\ref{24mar112}), we see that  in this region the function~$v(t,x)$ satisfies
\be\label{24mar304}
\bal
v_t-\kappa v_{xx}+\alpha_1(1-\bar z(t,x))v&=[\dot\xi(t)+\alpha_1(1-\bar z(t,x))\xi(t)]\un u(t,x)
\\
&\le [\dot\xi(t)+C\exp(-\omega t^\mu)\xi(t)]\un u(t,x)\le 0,
\enbal
\ee
if we choose $\xi(t)$ so that
\be
\dot\xi(t)+C\exp(-\omega t^\mu)\xi\le 0.
\ee
In particular, if $t_0>1$ is sufficiently large, we can take $\xi(t)$ as a solution to  
\be
\dot\xi(t)+\farc{m\xi(t)}{(1+t)^{1+m}}=0,
\ee
with some $m>0$ sufficiently small. This gives  
\be\label{24mar306}
\xi(t)= {\xi(t_0)}\exp\Big(\farc{1}{(1+t)^m}-\farc{1}{(1+t_0)^m}\Big).
\ee
Note that (\ref{24mar123}) holds with that choice.

{\bf Step 2. The intermediate region.} 
Next, we look at the intermediate region  $I_M$ where both terms in the right side of (\ref{24mar122}) are non-zero.
For the function $v_1(t,x)$ we have, as in (\ref{24mar304}):
\be\label{24mar308}
\bal
\pdr{v_1}{t}&-\kappa\pdrr{v_1}{x}+\alpha_1(1-\bar z(t,x))v_1=[\dot\xi(t)+\alpha_1(1-\bar z(t,x))]\un u(t,x)
\\
&\le [\dot\xi(t)+C\exp(-\omega t^\mu)\xi]\un u(t,x)\le -\farc{m\xi(t)}{2(1+t)^{1+m}}\un u(t,x)\le  -\farc{m\xi(0)}{4(1+t)^{1+m}}\un u(t,x).
\enbal
\ee
We used above the choice (\ref{24mar306}) for $\xi(t)$ as well as  (\ref{24mar123}).  Let us next compute
\be\label{24mar310}
\bal
\pdr{v_2}{t}-&\kappa\pdrr{v_2}{x}+\alpha_1[1-\bar z(t,x+c_*t)]v_2=
\farc{\beta-3/2}{t^{5/2-\beta}}
\cos\Big(\farc{x}{t^\mu}\Big)\\
&+\farc{\mu x}{t^{5/2-\beta+\mu}}\sin\Big(\farc{x}{t^\mu}\Big)
+\farc{\kappa}{t^{3/2-\beta+2\mu}}
\cos\Big(\farc{x}{t^\mu}\Big)+\alpha_1\farc{1-\bar z(t,x+c_*t)}{t^{3/2-\beta}}
\cos\Big(\farc{x}{t^\mu}\Big)\\
&=g_c(t,x)\cos\Big(\farc{x}{t^\mu}\Big)+g_s(t,x)\sin\Big(\farc{x}{t^\mu}\Big).
\enbal
\ee
Here, we have set
\be
\bal
&g_c(t,x)=-\farc{3/2-\beta}{t^{5/2-\beta}}+\farc{\kappa}{t^{3/2-\beta+2\mu}}+\alpha_1\farc{1-\bar z(t,x+c_*t)}{t^{3/2-\beta}},\\
&g_s(t,x)=\farc{\mu x}{t^{5/2-\beta+\mu}}.
\enbal
\ee
Thus, we have
\be\label{24jun1312}
\bal
{\cal L}[ v](t,x):&=\un v_t-\kappa\un v_{xx}+\alpha_1[1-\bar z(t,x+c_*t)]\un v\le  -\farc{m\xi(0)}{10(1+t)^{1+m}}\un u(t,x)
\\
&-g_c(t,x)\cos\Big(\farc{x}{t^\mu}\Big)
-g_s(t,x)\sin\Big(\farc{x}{t^\mu}\Big). 
\enbal
\ee
Let us consider separately each of the regions
\be
\bal
&I_{M,L}=\big\{L_1(t)\le x\le ({\pi}/{2})t^\mu\big\},~~I_{M,M}=\big\{({\pi}/2)t^\mu\le x\le {\pi}{t^\mu}\big\},\\
& I_{M,R}=\big\{{\pi}{t^\mu}\le x\le(3{\pi}/2){t^\mu}\big\}.
\enbal
\ee

In the region $I_{M,L}$, both the sine and the cosine appearing in the right side of (\ref{24mar310}) are positive. Moreover,
$g_s(t,x)>0$ as follows from its definition, and $g_c(t,x)>0$ as long as $t_0$ is sufficiently large and $\mu<1/2$. Here, we used 
(\ref{24mar112}) once again.  
Thus, we have 
\be
{\cal L}[v](t,x)\le 0,~~\hbox{ in the region $I_{M,L}$.} 
\ee

In the region $I_{M,M}$ the corresponding sine is still positive while cosine is negative. In addition, the last term
in the definition of $g_c(t,x)$ is controled by the exponentially decaying bound in (\ref{24mar112}). In that region, we also have a lower 
bound for $\un u(t,x)$  of the form
\be\label{24mar314}
\un u(t,x)\ge\farc{C_0x}{t^{3/2}}\ge \farc{C_0'}{t^{3/2-\mu}},
\ee
with the constants $C_0,C'$ that depend only on the initial condition $\phi(x)$ in (\ref{24mar302}). 
Therefore, we have
\be
\bal
{\cal L}[v](t,x):&\le  -\farc{m\xi(t_0)}{2(1+t)^{1+m}}\un u(t,x)
-\farc{2\kappa}{t^{3/2-\beta+2\mu}}\cos\Big(\farc{x}{t^\mu}\Big) \\
&\le
-\farc{C_0''\xi(t_0)}{(1+t)^{5/2+m-\mu}} 
+\farc{2\kappa}{t^{3/2-\beta+2\mu}}\le 0,~\hbox{ in the region $I_{M,M}$ for $t>t_0$.} 
\enbal
\ee
as long as $m>0$ and $\beta>0$ are sufficiently small, $t_0>1$ is sufficiently large, and $\mu>1/3$. 

In the region $I_{M,R}$ both the  sine and the  cosine in the right side of (\ref{24jun1312}) 
are negative, while~$\un u(t,x)$ still obeys the lower bound (\ref{24mar314}). 
Therefore, we have
\be
\bal
{\cal L}[v](t,x):&\le  -\farc{m\xi(t_0)}{2(1+t)^{1+m}}\un u(t,x)
-\farc{2\kappa}{t^{3/2-\beta+2\mu}}\cos\Big(\farc{x}{t^\mu}\Big)\\
&-\farc{\mu x}{t^{5/2-\beta+\mu}}\sin\Big(\farc{x}{t^\mu}\Big) \\
&\le
-\farc{C_0\xi(t_0)}{(1+t)^{5/2+m-\mu}} 
+\farc{2\kappa}{t^{3/2-\beta+2\mu}}+\farc{\mu }{t^{5/2-\beta }}\le 0,~\hbox{ in the region $I_{M,R}$,} 
\enbal
\ee
once again, as long as $m>0$ and $\beta>0$ are sufficiently small, $t_0>1$ is sufficiently large and $\mu>1/3$. 
 
It remains to check that the initial comparison (\ref{24mar318}) holds:
\be\label{24mar322}
P(t_0,x)=F(t_0,x)e^{\lambda_*(x-c_*t_0)}\ge  \max[0,\eps\un v(t_0,x)],~~\hbox{for all $x>0$.}
\ee
By choosing $\eps>0$ sufficiently small, we only need to ensure that (\ref{24mar322}) holds 
in the tails for $x\gg1$.  However, $P(t_0,x)$ satisfies a lower bound
\be
P(t_0,x)\ge e^{-\alpha_1t_0}\bar u(t,x).
\ee
Here, $\bar u(t,x)$ is the solution to the heat equation on the whole line:
\be
\bal
&\bar u_t=\kappa\bar u_{xx},
&\bar u(0,x)=P(0,x).
\enbal
\ee
If we take $\un u(0,x)=P_0(x)$, we will ensure that $\bar u(t,x)\ge \un u(t,x)$ for all $t>0$ and $x>0$. Thus, we can take
$\phi(x)=P(0,x)$ in (\ref{24mar302}) and then choose $\eps>0$ sufficiently small to guarantee that 
(\ref{24mar322}) holds. 

Summarizing, we have shown, in particular, that $P(t,x)$ satisfies
a lower bound
\be
P(t,x)\ge \eps\un u(t,x),~~x>\sqrt{t},~~t>t_0.
\ee
This, in turn, implies that
\be
\label{eq:rec}
Q(t,x)\ge \un Q(t,x):=\eps\un u(t,x-c_*t)e^{-\lambda_*(x-c_*t)},~~x>c_*t+\sqrt{t}.
\ee
 
\subsubsection{Asymptotics of the learning front for the sub-solution}

It remains to show that the learning front of the sub-solution $\un Q(t,x)$ satisfies
\be\label{24mar402}
\tel (t;\un Q)\ge (c_*+\delta)t-L_1,~~\hbox{ for all $t>t_0$},
\ee
and
also that 
\be\label{24mar414}
\tel(t;\un Q)\ge (\kappa+\alpha_1-\delta)t+o(t).
\ee
To this end, recall~\eqref{eq:rec} 
so that the function $J(t,x)$ defined by~\eqref{24feb704} satisfies 
\be
\bal
J(t,c_*t+x)&\ge  J[\un Q](t, c_*t+x):=\farc{\eps e^{-(x+c_*t)}}{\rho-\kappa}\int_{x+c_*t}^\infty e^y\un u(t,y-c_*t)
e^{-\lambda_*(y-c_*t)}dy\\
&=\farc{\eps e^{-x}}{\rho-\kappa}\int_x^\infty e^{(1-\lambda_*)y}\un u(t,y)dy.
\enbal
\ee 
The function  $\un u(t,x)$ has the form
\be
\bal
\un u(t,x)&=\farc{1}{\sqrt{4\pi t}}\int_0^\infty
\Big(e^{-(x-y)^2/(4\kappa t)}-e^{-(x+y)^2/(4\kappa t)}\Big)\phi(y)dy\\
&=
\farc{e^{-|x|^2/(4\kappa t)}}{\sqrt{4\pi \kappa t}}
\int_0^\infty e^{-|y|^2/(4\kappa t)}
\Big(e^{xy/(2\kappa t)}-e^{-xy/(2\kappa t)}\Big)\phi(y)dy.
\enbal
\ee
We may assume without loss of generality that the initial condition $\phi(y)$ is supported inside the
interval $1\le y\le 2$. In that case, if $x>c_*t$, we have the inequality
\be
e^{xy/(2\kappa t)}\ge\delta_1 e^{-xy/(2\kappa t)}.
\ee
We deduce that for $x>c_*t$ we have 
\be
\bal
\un u(t,x)&\ge \farc{\delta_2}{\sqrt{t}}\int_0^\infty
e^{-(x-y)^2/(4\kappa t)}\phi(y)dy \ge \farc{\delta_3}{\sqrt{t}}e^{-x^2/(4\kappa t)}.
\enbal
\ee
It follows that for all $x>0$ we have
\be
\bal
J(t,c_*t+x)&\ge \farc{e^{-x}}{C\sqrt{t}} \int_x^\infty e^{(1-\lambda_*)y}e^{-y^2/(4\kappa t)} dy
= \frac{e^{-x+(1-\lambda_*)^2\kappa t}}{C\sqrt{t}}
\int_{x}^\infty e^{-|y/(2\sqrt{\kappa t})-(1-\lambda_*)\sqrt{\kappa t}|^2 }dy\\
&= Ce^{-x+(1-\lambda_*)^2\kappa t}\int_{\Gamma(t)}^\infty e^{-|y|^2 }dy,
\enbal
\ee
with
\be
\Gamma(t)=\farc{x}{\sqrt{4\kappa t}}-(1-\lambda_*)\sqrt{\kappa t}.
\ee
Therefore, we have $J(t,x+c_*t)\gg 1$ if we have both 
\be\label{24mar410}
-x+(1-\lambda_*)^2\kappa t\gg 1,
\ee
and $\Gamma(t)\le 0$, which is
\be\label{24mar412}
x<2(1-\lambda_*)\kappa t
\ee
Note that (\ref{24mar412}) is a consequence of (\ref{24mar410}) when $\lambda_*<1$. 
In particular, we have 
\be
J(t,x+(c_*+\delta)t))\gg \farc{1}{\alpha'(1)}=\farc{2}{\alpha_1},
\ee
as long as 
\be
\delta<\kappa(1-\lambda_*)^2. 
\ee
On the other hand, we have
\be
c_*+(1-\lambda_*)^2\kappa =2\sqrt{\kappa\alpha_1}+(\sqrt{\kappa}-\sqrt{\alpha_1})^2=\kappa+\alpha_1.
\ee
This implies both (\ref{24mar402}) and (\ref{24mar414}) finishing the proof of Lemma~\ref{lem-mar304}.~$\Box$ 

\begin{appendix}

\section{The derivative formulation of the Lucas-Moll system}\label{sec:append-deriv} 

Here, we derive the reformulation (\ref{eqn:intro1}) of the Lucas-Moll system 
as a system of two semilinear parabolic equations, one being forward in time, and
the other  backward in time. This reformulation is motivated by the
construction of the balanced growth paths in~\cite{Porretta-Rossi}. 

We start with the  Lucas-Moll system written in the form  (\ref{Psi2})-(\ref{V2}):
\begin{equation}\label{23dec2702}
\begin{aligned}
\pdr{F(t,x)}{t}&=\kappa \pdrr{F(t,x)}{x}+F(t,x)\int_{-\infty}^x \alpha(s^*(t,y))(-F_y(t,y))dy,\\
\rho V(t,x)&=\pdr{V(t,x)}{t}+\kappa\pdrr{V(t,x)}{x}
+\max_{s\in[0,1]}\Big[(1-s)e^x+\alpha(s)\int_x^\infty [V(t,y)-V(t,x)] (-F_y)(t,y)dy\Big].
\end{aligned}
\end{equation}
The optimal control $s^*(t,x)$ that appears in the first equation in (\ref{23dec2702}) is
\be\label{23dec2704}
\bal
s^*(t,x)&=\argmax_{s\in[0,1]}
\Big[(1-s)e^x+\alpha(s)\int_x^\infty [V(t,y)-V(t,x)] (-F_y)(t,y)dy\Big]
\\
&=\argmax_{s\in[0,1]}
\Big[(1-s) +\alpha(s)I(t,x)\Big].
\enbal
\ee 
Here, we have defined the integral that appears in the definition of $s^*(t,x)$ as 
\be\label{24jan2912}
I(t,x):=e^{-x}\int_x^\infty [V(t,y)-V(t,x)] (-F_y)(t,y)dy=e^{-x}\int_x^\infty v(t,y)F(t,y)dy,
\ee
where
\be\label{24apr2602}
v(t,x)=V_x(t,x).
\ee

In order to describe the maximizer in (\ref{23dec2704}), 
let us for the moment assume that $V(t,x)$ is increasing so that
$v(t,x)>0$ -- we will come back to this point later. 
Then,~$I(t,x)$ is a positive strictly decreasing function such that 
\be
\bal
& I(t,x)\to 0~~\hbox{ 
as~$x\to+\infty$},\\
& I(t,x)\to +\infty \hbox{ as~$x\to-\infty$. }
\enbal
\ee
Thus, there exists a 
unique point $\eta_\ell(t)\in\Rm$ so that 
\be\label{23dec2721}
I(t,\eta_\ell(t)) =\farc{1}{\alpha'(1)},
\ee
and
\be
\bal
&I(t,x)<\farc{1}{\alpha'(1)},~~\hbox{ for $x>\eta_\ell(t)$},\\
&I(t,x)>\farc{1}{\alpha'(1)},~~\hbox{ for $x<\eta_\ell(t)$.}
\enbal
\ee
Let us also define the function $s_m(I)$ as follows: first, we set 
\be\label{24jan3102}
s_m(I)=1,~~\hbox{ for }~~I>\farc{1}{\alpha'(1)}.
\ee
Second, for $I\le 1/\alpha'(1)$ the function $s_m(I)$ is determined implicitly by
\be\label{24jan2908}
\alpha'(s_m(I))=\farc{1}{I},~~\hbox{for}~~I\le \farc{1}{\alpha'(1)}.
\ee
With this definition, the maximizer in (\ref{23dec2704}) can be computed explicitly 
and is given by 
\be\label{24jan2904}
s^*(t,x)=s_m(I(t,x)).
\ee

As, by assumption, $\alpha(s)$ is concave, the function $\alpha'(s)$ is decreasing in $s\ge 0$.
Thus, the function~$s_m(I)$ is increasing in $I$. Since~$I(t,x)$ is decreasing in $x$, 
it follows that $s^*(t,x)$ is decreasing in~$x$, as expected. 
In particular, we have 
\be\label{24feb502}
\hbox{$s^*(t,x)=1$ for all $x<\eta_\ell(t)$ and $s^*(t,x)<1$ 
for all $x>\eta_\ell(t)$.}
\ee
Because of that property, we will refer to $\eta_\ell(t)$, determined by~(\ref{23dec2721}), 
as the learning
front location.   

We now obtain an equation for the derivative $v(t,x)$ defined in (\ref{24apr2602}). The above notation allows us to write the 
equation for the function $V(t,x)$ in (\ref{23dec2702}) as
\be\label{24jan2910}
\bal
\rho V(t,x)&=\pdr{V(t,x)}{t}+\kappa\pdrr{V(t,x)}{x}+e^x\Big[1-s_m(I(t,x))+\alpha(s_m(I(t,x)))I(t,x)\Big].
\enbal
\ee
Differentiating (\ref{24jan2910}) in $x$ gives
\be\label{24jan2914}
\bal
\rho v(t,x)&=\pdr{v(t,x)}{t}+\kappa\pdrr{v(t,x)}{x}+e^x\Big[1-s_m(I(t,x))+\alpha(s_m(I(t,x)))I(t,x)\Big]\\
&+e^x\Big[-s_m'(I(t,x))I_x(t,x)+\alpha'(s_m(I(t,x)))s_m'(I(t,x))I_x(t,x)I(t,x)+\alpha(s_m(I(t,x)))I_x(t,x)\Big].
\enbal
\ee
Recall that, by (\ref{24jan2908}), we have
\be
\alpha'(s_m(I(t,x)))I(t,x)=1,~~\hbox{ for $x>\eta_\ell(t)$}.
\ee
Furthermore, by the definition (\ref{24jan2912}) of $I(t,x)$, we have
\be\label{24apr2606}
I_x(t,x)=-I(t,x)-e^{-x}v(t,x)F(t,x).
\ee
Thus, for $x>\etal(t)$ we can write the term inside the parentheses in the last line of (\ref{24jan2914}) as
\be
\label{23dec2726}
\bal
&-s_m'(I(t,x))I_x(t,x)+\alpha'(s_m(I(t,x)))s_m'(I(t,x))I_x(t,x)I(t,x)+\alpha(s_m(I(t,x)))I_x(t,x)\\
&=\alpha(s_m(I(t,x)))I_x(t,x)
=-\alpha(s_m(I(t,x)))I(t,x)-\alpha(s_m(I(t,x)))e^{-x}v(t,x)F(t,x).
\enbal
\ee
Using this in (\ref{24jan2914}) gives
\be\label{23dec2727}
\bal
\rho v(t,x)&=\pdr{v(t,x)}{t}+\kappa\pdrr{v(t,x)}{x}+e^x\Big[1-s_m(I(t,x))+\alpha(s_m(I(t,x)))I(t,x)\Big]\\
&+e^x\Big[-\alpha(s_m(I(t,x)))I(t,x)-\alpha(s_m(I(t,x)))e^{-x}v(t,x)F(t,x)\Big]\\
&=\pdr{v(t,x)}{t}+\kappa\pdrr{v(t,x)}{x}+e^x\Big[1-s_m(I(t,x)) -\alpha(s_m(I(t,x)))e^{-x}v(t,x)F(t,x)\Big], 
\enbal
\ee
for $x>\etal(t)$.
%

On the other hand, for $x<\etal(t)$, we have 
\be\label{24apr2608}
s_m(I(t,x))=1,
\ee
and  in that region  (\ref{24jan2914}) takes the form
\be\label{24apr2604}
\bal
\rho v(t,x)&=\pdr{v(t,x)}{t}+\kappa\pdrr{v(t,x)}{x}+e^x\alpha(1)I(t,x) \
+e^x\alpha(1)I_x(t,x) \\
&=\pdr{v(t,x)}{t}+\kappa\pdrr{v(t,x)}{x}-\alpha(1)v(t,x)F(t,x),
\enbal
\ee
which is identical to (\ref{23dec2727}) because of (\ref{24apr2608}). Thus, $v(t,x)$
satisfies (\ref{23dec2727}) everywhere.

Next, it is helpful to take out the exponential factor and represent $v(t,x)$ as 
\be\label{24feb516}
v(t,x)=\farc{1}{\rho-\kappa}w(t,x)e^x.
\ee
Inserting this form into (\ref{23dec2727}), we see that the function $w(t,x)$ satisfies
\be\label{24jan3002}
\bal
\rho w(t,x)
&=\pdr{w(t,x)}{t}+\kappa\pdrr{w(t,x)}{x}+2\kappa w_x(t,x)+\kappa w(t, x)+(\rho-\kappa)(1-s_m(I(t,x)))\\
& -\alpha(s_m(I(t,x)))w(t,x)F(t,x).
\enbal
\ee
This equation holds for $x>\etal(t)$.

With these changes of variables, the Lucas-Moll system takes the form of a coupled system of forward and backward
nonlocal semilinear parabolic equations
\be\label{24jan3004} 
\bal
&\pdr{F(t,x)}{t}=\kappa \pdrr{F(t,x)}{x}+F(t,x)\int_{-\infty}^x\alpha(s_m(I(t,y))(-F_y)(t,y)dy,\\
&\pdr{w(t,x)}{t}+\kappa\pdrr{w(t,x)}{x}+2\kappa w_x(t,x)+(\rho-\kappa)\big[1-s_m(I(t,x))-w(t,x)\big]\\
&~~~~~~~~~~~-\alpha(s_m(I(t,x)))w(t,x)F(t,x)=0,
\\
&I(t,x)=\farc{1}{\rho-\kappa}e^{-x}\int_x^\infty e^yw(t,y)F(t,y)dy.
\enbal
\ee
The function $F(t,x)$ satisfies the boundary conditions
\be\label{24jan3104}
F(t,-\infty)=1,~~F(t,+\infty)=0.
\ee
while $w(t,x)$ satisfies the boundary conditions
\be\label{24jan3008}
w(t,-\infty)=0,~~w(t,+\infty)=1.
\ee
The system  (\ref{24jan3004})  has the form (\ref{eqn:intro1})

\section{The form of the equations for $\alpha(s)=\alpha_1s^k$ }\label{sec:alpha-powers}

To simplify slightly some considerations, we often  assume that  $\alpha(s)$ has the form
\be\label{24jan3010}
\alpha(s)=\alpha_1s^k,
\ee
with some $k\in[1/2,1)$. 
Then, $s_m(I)$ is determined by (\ref{24jan2908}) that takes the form   
\be
\farc{k\alpha_1}{s_m^{1-k}(I)}= {I}^{-1},~~\hbox{for } I\le\farc{1}{k\alpha_1},
\ee
so that 
\be\label{24jan3012}
\bal
&s_m(I)=(k\alpha_1I)^{1/(1-k)},
~~\hbox{for } I\le\farc{1}{k\alpha_1},\\
&s_m(I)=1,~~\hbox{for } I\ge\farc{1}{k\alpha_1},
\enbal
\ee
and
\be\label{24apr2618}
\bal
&\alpha(s_m(I))=\alpha_1(k\alpha_1I)^{k/(1-k)},~~\hbox{for } I\le\farc{1}{k\alpha_1},\\
&\alpha(s_m(I))=\alpha_1,~~\hbox{for } I\ge\farc{1}{k\alpha_1},
\enbal
\ee
The transition point $\eta_\ell(t)$ is determined in this case by the condition
\be\label{24jan3014}
I(t,\eta_\ell(t))=\farc{1}{k\alpha_1}.
\ee
In particular, when $\alpha(s)=\alpha_1\sqrt{s}$, we have 
\be\label{24apr1202}
\bal
&s_m(I)=\Big(\farc{\alpha_1I}{2}\Big)^{2},
~~\hbox{for } I\le\farc{2}{\alpha_1},\\
&s_m(I)=1,~~\hbox{for } I\ge\farc{2}{\alpha_1},
\enbal
\ee
and
\be\label{24apr1204}
\bal
&\alpha(s_m(I))=\farc{\alpha_1^2I}{2},~~\hbox{for } I\le\farc{2}{\alpha_1},\\
&\alpha(s_m(I))=\alpha_1,~~\hbox{for } I\ge\farc{2}{\alpha_1},
\enbal
\ee
The transition point $\eta_\ell(t)$ is determined in this case by the condition
\be\label{24apr1206}
I(t,\eta_\ell(t))=\farc{2}{\alpha_1}.
\ee

%

\end{appendix}

\bibliographystyle{plain}

\end{document}